\documentclass[reqno,11pt]{amsart}
\usepackage{amsfonts,amsmath,amsthm, url}
\usepackage{amssymb}  
\usepackage{mathrsfs}
\usepackage{graphics,color}
\usepackage{graphicx}
\usepackage{stmaryrd} 
\usepackage{mathrsfs}
\usepackage[abbrev]{amsrefs}
\usepackage{comment}

   \newtheorem{thm}{Theorem}[section]
   \newtheorem{prop}{Proposition}[section]
   \newtheorem{lem}{Lemma}[section]

   \newtheorem{defn}{Definition}[section]

\newcommand{\N}{\mathbb{N}}
\newcommand{\Z}{\mathbb{Z}}

\newcommand{\R}{\mathbb{R}}

\newcommand{\ol}{\varolessthan}
\newcommand{\oo}{\varodot}
\newcommand{\og}{\varogreaterthan}

\makeatletter
  
  \@addtoreset{equation}{section}
\makeatother

\newcommand{\black}{\color{black}}

\begin{document}

\title{
Besov spaces associated with the Harmonic oscillator}
\author[Reika Fukuizumi]{Reika Fukuizumi$^{\scriptsize 1}$} 
\author[Tsukasa Iwabuchi]{Tsukasa Iwabuchi$^{\scriptsize 2}$}

\date{\today}

\keywords{harmonic oscillator, Besov space}
\subjclass{30H25}

\maketitle

\begin{center} 
$^1$ Department of Mathematics, Faculty of Fundamental Science and Engineering, Waseda University,\\
169-8555 Tokyo, Japan; \\
\email{fukuizumi@waseda.jp}
\end{center}

\begin{center} 
$^2$ Mathematical Institute, Faculty of Science, Tohoku University,\\
980-8578, Sendai, Japan; \\
\email{t-iwabuchi@tohoku.ac.jp}
\end{center}

\begin{abstract}
The Besov space associated with the harmonic oscillator is introduced and thoroughly explored in this paper. It provides a comprehensive summary of the fundamental concepts of the Besov spaces, their embedding properties, bilinear estimates, and related topics.
\end{abstract}
\black

\section{Introduction}
We study the Besov space based on the Littlewood-Paley decomposition
associated with the harmonic oscillator on $\mathbb R^d$, for $d \geq 1$, 
\[
H = -\Delta + |x|^2 . 
\]
The operator $H$ is one of the important operators in quantum mechanics. Moreover, \black when rigorously analyzing physically significant nonlinear equations, for example, the Gross-Pitaevskii equation \cites{dbdf1, dbdf2, dbdf3, t, Y-2004}, the Sobolev spaces and Besov spaces based on this \black harmonic oscillator as the fundamental operator, as well as the bilinear estimates in these spaces, are extremely useful.

The eigenvalues of $H$ are well known, and the eigenfunctions are written explicitly using Hermite functions. In this paper, we decompose the spectrum of $H$ to introduce dyadic decomposition, 
and utilize the boundedness of the spectral multiplier to introduce Besov spaces associated with the operator $H$. 
The aim of this paper is to establish basic estimates in the Besov spaces associated with the operator $H$. 

The Hermite Besov spaces have been introduced by Petrushev and Xu~\cite{PX-2008} (see also \cites{BD-2015,BDY-2012}), in a different way from this paper, 
 based on the 
Calder\'on reproducing formula for the identity operator. 
The spaces introduced by them are equivalent to ours  
(see Theorem \ref{thm:equiv} below). 
\black 
Since we prefer the setting better adapted to the analysis of partial differential equations, in this paper, we introduce Besov spaces associated with $H$ following the argument in \cite{IMT-2019}, whose key feature is that it deals with Besov spaces based on the Dirichlet Laplacian. 

We see that $H$ has a self-adjoint realization on $L^2 $, 
and can be written as follows. 
\[
\begin{cases}
D(H) = \{ f \in L^2 \, | \, \Delta f , |x|^2 f \in L^2  \}, 
\\ 
Hf = (-\Delta + |x|^2) f, \quad f \in D(H). 
\end{cases}
\]
By applying the spectral theorem, 
the resolution of the identity $\{ E (\lambda) \}_{\lambda \in \mathbb R}$ 
exists such that 
\[
\begin{split}
& \lim_{\lambda \to \infty} E(\lambda) f  = f \text{ in } L^2, 
\quad \text{for all } f \in L^2 ,
\\
& Hf = \int_{-\infty}^\infty \lambda \, dE(\lambda)f 
\text{ in } L^2 , \quad f \in D(H). 
\end{split}
\]
We denote $\sigma(H)$ the spectrum of $H$. It is known that 
$\inf \sigma(H)$ is strictly positive, \black which implies the equivalence between the two norms 
of the homogeneous and the non-homogeneous types, 
\[
\| H f \|_{L^2} \simeq \| (1+H)f \|_{L^2}.
\] 
Remark that 
\[
\| f \|_{L^2} 
\leq 
\frac{1}{ \inf \sigma (H) } \| Hf \|_{L^2}. 
\]
Because of this, the two norms define one function space, 
and we will use the left hand side to introduce Besov spaces. 
We take $\phi_0 \in C_0^\infty (\mathbb R)$ a non-negative function on 
$\mathbb R$ such that 
\[
{\rm supp \, } \phi_0 \subset [2^{-1},2], \quad 
\sum _{ j \in \mathbb Z} \phi_0 (2^{-j} \lambda) = 1 
 \text{ for } \lambda > 0 ,
\]
and $\{ \phi_j \}_{j \in \mathbb Z}$ is defined by 
\[
\phi_j(\lambda) := \phi_0 (2^{-j}\lambda) , \quad \text{ for } 
\lambda \in \mathbb R . 
\]

\vskip3mm 




\vskip3mm 

\begin{defn} For $s \in \mathbb R$, 
$1 \leq p,q \leq \infty$, $B^s_{p,q}(H)$ is defined by 
\[
B^s_{p,q}(H) := 
\{ f \in \mathcal S'(\mathbb R^d) \, | \, \| f \|_{B^s_{p,q}(H)} < \infty 
\},
\]
where 
\[
\| f \|_{B^s_{p,q}(H)} := 
\Big\| \Big\{ 2^{sj} \| \phi_j(\sqrt{H}) f \|_{L^p}
\Big\}_{j=-1}^\infty
\Big\|_{\ell^q } . 
\]
\end{defn}
\vskip3mm 

The first number $j=-1$ in the sequence is determined by $j_0 \in \mathbb Z$ such that

\begin{equation}\label{0321-2}
2^{j_0+1} \le \inf \sigma (H). 
\end{equation}
 \black

For simplicity, we will write the sum over 
$j \in \mathbb Z$, and explicitly indicate the sum over 
$j$ with $j \geq -1$ when a clarification is needed. 
On the partition of the unity it reads that 
\[
f = \sum _{ j \in \mathbb Z} \phi_j (\sqrt{H}) f 
= \sum _{j = -1}^\infty \phi_j (\sqrt{H}) f, 
\]
since 
\[
\phi_j (\sqrt{H})f = 0  \text{ if  }j < -1. 
\]
\black 
We notice that the positivity of the spectrum of $H$ implies 
the following equivalence. 
\[
\Big\| \Big\{ 2^{sj} \| \phi_j(\sqrt{H}) f \|_{L^p}
\Big\}_{j=-1}^\infty
\Big\|_{\ell^q }
\simeq 
\| \psi(\sqrt{H})f \|_{L^p} + 
\Big\| \Big\{ 2^{sj} \| \phi_j(\sqrt{H}) f \|_{L^p}
\Big\}_{j \in \mathbb N}
\Big\|_{\ell^q },
\]
where 
$\psi \in C_0^\infty (\mathbb R)$ satisfies 
$\psi + \sum_{j=1}^\infty \phi_j = 1$. 

\vskip3mm

Let us introduce the basic properties of the Besov space 
$B^s_{p,q}(H)$ in the following proposition.  

\vskip3mm 

\begin{prop}\label{prop:intro1}
Let $s, \alpha \in \mathbb R$ and $1 \leq p,q,r \leq \infty$. 
The following {\rm (i)-(vii)} hold: 
\begin{enumerate}

\item[(i)] $B^s_{p,q}(H)$ is a Banach space and enjoys 
$\mathcal S(\mathbb R^d) \hookrightarrow B^s_{p,q}(H) \hookrightarrow \mathcal S'(\mathbb R^d)$. 

\item[(ii)] If $1\le p,q < \infty$ \black and 
$1/p + 1/p' = 1/q + 1/q' = 1$, then the dual space of 
$B^s_{p,q}(H)$ is $B^{-s}_{p',q'}(H)$. 
Moreover, for any $p,q \in [1,\infty]$, 
we have the following norm equivalence. 
\[
\| f \|_{B^{-s}_{p',q'}}
\simeq 
\sup_{\| g \|_{B^s_{p,q}}=1} 
\Big| 
\sum _{j=-1}^\infty
\int_{\mathbb R^d} \phi_j(\sqrt{H})f(x) \overline{\Phi_j(\sqrt{H})g(x)} ~dx
\Big|,
\]
where $\Phi_j = \phi_{j-1} + \phi_j + \phi_{j+1} $. 
Denote $Q^{-s}_{p',q'}:=\{f \in \mathcal{S}, \|f\|_{B^{-s}_{p',q'}} \le 1\}$. 
If $g \in \mathcal{S}'$, we then have 
$$ \black \|g\|_{B^s_{p,q}} \le C \sup_{f \in Q_{p',q'}^{-s}} |\langle f,g\rangle |.$$


\item[(iii)] If $r \leq p$, then 
$B^{s+d(\frac{1}{r}-\frac{1}{p})}_{r,q}(H) 
\hookrightarrow B^s_{p,q}(H)$. 

\item[(iv)] For every $f \in B^{s+\alpha}_{p,q}(H)$, 
$H^{\frac{\alpha}{2}}f \in B^s_{p,q}(H)$.

\item[(v)] If $s<\alpha$, the space  $B^{\alpha}_{p,\infty}(H)$ is compactly embedded into $B^{s}_{p,1}(H).$

\item[(vi)] There exists a constant $C>0$ such that 
$$ C^{-1} \|f\|_{B^0_{p,\infty}(H)} \le \|f\|_{L^p} \le C\|f\|_{B^0_{p,1}(H)}$$

\item[(vii)] 
Let $s, s_0 >0$, $p, r, r_0 \in [1, \infty]$, $\theta \in(0,1)$ satisfy 
\begin{equation}\notag 
s- \frac{d}{p} 
= \theta \Big( - \frac{d}{r} \Big) + (1-\theta) \Big( s_0 - \frac{d}{r_0} \Big) ,
\end{equation}
\begin{equation}\notag 
-\frac{d}{r} \not = s_0 - \frac{d}{r_0}, 
\qquad 
\begin{cases}
s \leq (1-\theta) s_0 
& \text{if } \max\{r, r_0\} \leq p,
\\
s < (1-\theta) s_0
& \text{if } \min\{r, r_0 \} \leq p  < \max\{r, r_0\}. 
\end{cases}
\end{equation}
Then we have 
$$\|f\|_{B^{s}_{p,1}(H)} \le \|f\|^{\theta}_{B^{0}_{r,\infty}(H)}
\|f\|^{1-\theta}_{B^{s_0}_{r_0,\infty}(H)}.
$$
\end{enumerate}
\end{prop}
\vskip3mm 

Remark that the above items (i)--(iv), 
where $H$ is replaced by the Dirichlet Laplacian, 
have been already established in \cite{IMT-2019}, and those arguments can be applied similarly for the case $H$. The equivalent norm in (ii)  follows from the property of duality. 
We will thus give a brief proof only for (v), (vi) and (vii) in this paper, in Appendix. \black 
\vskip3mm

We here mention that the uniform boundedness of the operators
$\phi_j(\sqrt{H})$ in $j$ holds.
\begin{equation}\label{unifbdd}
\begin{split}
\sup_{j \geq -1} \| \phi (2^{-j}\sqrt{H}) \|_{L^p \to L^p} 
< 
\infty , 
\end{split}
\end{equation}

for all $1 \leq p \leq \infty $ and 
$ \phi \in \mathcal S(\mathbb R)$. This holds by the same reason 
as in Section 8 in~\cite{IMT-RMI}. 
Initially, $\phi_j(\sqrt{H})$ 
is defined on $L^2$ with an application of the spectral theorem
and is a bounded operator on $L^2$. This uniform boundedness (\ref{unifbdd}) on $L^p$, 
$1 \leq p \leq \infty$ plays a very important role
to establish the theory of Besov spaces, 
we thus give a brief proof in the appendix. 
\black 

\vskip3mm 

We write the Bony paraproduct formula. 
\[
\begin{split}
fg =& \sum_{k} \phi_k(\sqrt{H})f 
\sum_{l } \phi_l (\sqrt{H}) g
= 
 \Big( \sum_{k \leq l-2} + \sum _{l \leq k-2} + \sum _{|k-l| \leq 1} \Big) 
  ( f_k g_l ),
\\
=& f \ol g + f \og g + f \oo g . 
\end{split}
\]
where $f_k = \phi_k(\sqrt{H})f $ and $g_l = \phi_l (\sqrt{H}) g$. 
Then, we have the bilinear estimates as follows. 
\vspace{3mm}

\begin{prop}\label{prop:bilinear} Let $s, r \in \mathbb R, 1 \leq p,p_1,p_2,q \leq \infty$ and $1/p = 1/p_1 + 1/p_2$. 
\begin{itemize}
\item[(i)]  If $s>0$, then \black
there exists a constant $C>0$ such that 
\begin{eqnarray}\label{est:lowhigh1}
\| f\ol g \|_{B^s_{p,q}(H)} 
\leq C \| f \|_{L^{p_1}} \| g \|_{B^s_{p_2,q}(H)}.
\end{eqnarray}
\item[(ii)] If $s<0$  and $s+r>0$, then \black   
$$ \| f\ol g \|_{ B^{s+r}_{p,q}(H)} 
\leq C \| f \|_{ B^{s}_{p_1, \infty} } \| g \|_{ B^r_{p_2,q}(H)}.$$
\item[(iii)] If $s = s_1 + s_2 > 0$, $1\leq q_1, q_2 \leq \infty$ and 
$1/q = 1/q_1 + 1/q_2$, then 
\begin{eqnarray} \label{est:resonant}
\| f\oo g \|_{ B^s_{p,q}(H)} 
\leq C \| f \|_{ B^{s_1}_{p_1,q_1}(H)} \| g \|_{B^{s_2}_{p_2,q_2}(H)}.
\end{eqnarray}
\end{itemize}
\end{prop}

\noindent {\bf Remark.} In the definition of the para product above we
divided into the cases: 
$${k \leq l-2}, \quad {l \leq k-2}, \quad {|k-l| \leq 1}$$ 
but any number $N_0 \in \N$ for this division works for the proof, for example we may consider 
$${k < l-N_0}, \quad {l > k-N_0}, \quad {|k-l| \leq N_0}.$$
\black

\noindent{\bf Remark.} 
The regularity assumptions $s>0$ in (i) and $s+r>0$ in (ii) are likely to be technical rather than intrinsic, 
and therefore may not be optimal, since no such restrictions are needed in the classical Laplacian case.
\black
\vskip3mm

As a simple application of this Proposition \ref{prop:bilinear}, 
we have the following bilinear estimates.
\vspace{3mm}

\begin{thm} \label{thm:bilinear}
\begin{itemize}
\item[(i)]
Let  
\[
s > 0 , \quad 1 \leq p, p_1, p_2, p_3, p_4 , q\leq \infty , 
\quad \dfrac{1}{p}= \dfrac{1}{p_1} + \dfrac{1}{p_2}
= \dfrac{1}{p_3} + \dfrac{1}{p_4}.
\] 
Then there exists a positive constant $C$ such that
\[
\| fg \|_{B^s_{p,q}(H)} 
\leq C 
\Big( \| f \|_{ B^s_{p_1,q}(H)} \| g \|_{L^{P_2}}
    + \| f \|_{L^{p_3}} \| g \|_{ B^s_{p_4,q}(H)} 
\Big) ,
\]
for all $f \in B^s_{p_1,q}(H) \cap L^{p_3}, 
g \in L^{p_2} \cap B^s_{p_4,q}(H)$. 
\item[(ii)] Let $s<0<r$, $s+r>0$ and 
$1 \leq p,p_1,p_2\leq \infty$ with $\frac{1}{p}=\frac{1}{p_1}+\frac{1}{p_2}.$ Then, we have 
$$ \|fg\|_{B^s_{p,q}} \le C \|f\|_{B^{s}_{p_1,q}} \| g\|_{B^{r}_{p_2,q}},$$
for $f\in B^{s}_{p_1,q}(H)$ and $g\in B^{r}_{p_2,q}(H)$. 
\end{itemize}
\end{thm}
\vspace{3mm}

\noindent {\bf Remark.} Indeed, by Proposition \ref{prop:bilinear}, 
we can estimate each paraproduct under the parameters' condition of (ii) as follows. 
\begin{eqnarray*}
&&\|f\ol g\|_{B^{s}_{p,q}} \lesssim 
\|f\ol g\|_{B^{s+r}_{p,q}} 
\lesssim \|f\|_{B^{s}_{p_1,q}}\|g\|_{B^{r}_{p_2,q}}, \\
&& \|f\oo g\|_{B^{s}_{p,q}} \lesssim
\|f\oo g\|_{B^{s+r}_{p,q}} 
\lesssim \|f\|_{B^{s}_{p_1,q}}\|g\|_{B^{r}_{p_2,q}}, \\
&&\|f \og g\|_{B^{s}_{p,q}} 
\lesssim \|f\|_{B^{s}_{p_1,q}}\|g\|_{B^{r}_{p_2,q}}. 
\end{eqnarray*}
\black

\vskip2mm 

\noindent {\bf Remark.} Only for the purpose to give a proof for Theorem \ref{thm:bilinear}, it is sufficient to use 
the decomposition $fg$ into two parts (see the proof):
$$ fg = \sum_{k} \phi_k(\sqrt{H})f 
\sum_{l } \phi_l (\sqrt{H}) g
=
 \Big( \sum_{k \leq l} + \sum _{l < k} \Big) ( f_k g_l ).$$
\black 

\vskip3mm

\noindent {\bf Remark. } The inequality {\rm (i)} 
for the Besov spaces associated with the Laplacian is well-known (see e.g.~\cite{RS_1996}). 
In the Sobolev spaces associated to $H$, the existing estimate is as follows. 
\[
\| H^s (fg) \|_{L^p} 
\leq C \Big( \| H^s f \|_{L^{p_1}} \| g \|_{L^{p_2}} 
   + \| f \|_{L^{p_3}} \| H^s g \|_{L^{p_4}}
\Big) ,
\]
where $s > 0$, $1 < p , p_1,p_2, p_3, p_4 < \infty$ 
and $1/p = 1/p_1 + 1/p_2 = 1/p_3 + 1/p_4$ 
(see~\cite{t}). 
This is proved by the following equivalence
between the norms (\cite{dg}, also see Proposition \ref{prop:1}), 
\[
\| H^s f \|_{L^p} 
\simeq \| |x|^{2s} f \|_{L^p} + \| (-\Delta)^{s} \|_{L^p}, 
\quad s> 0 , 1 < p < \infty,
\]
using the H\"older inequality and 
the bilinear estimate for the standard Laplacian $(-\Delta)^s$ 
(see e.g., ~\cite{taylor}). 
We underline that in the Besov spaces $B^s_{p,q}(H)$, it is possible to include the {\it indices 
$p = 1$ and $\infty$}, and the present paper gives a proof for this fact. 


\vskip3mm 

Following the similar arguments in the paper~\cite{Iw-2018}, 
we have the following results about the smoothing effects of the semigroup $\{ e^{-tH} \}_{t\geq 0}$. 

\begin{thm}\label{thm:bdd}
Let $t \geq 0$, 
$s, s_1 ,s_2 \in \mathbb R$, $1 \leq p,p_1,p_2, q, q_1 , q_2 \leq \infty$. 
\\
{\rm (i)} 
$e^{-tH}$ is a bounded linear operator in 
$B^s_{p,q } (H)$, i.e., 
there exists a constant $C > 0$ such that 
 for any $f \in B^s_{p,q} (H)$
\begin{equation}\label{317-11}
 e^{-tH} f \in B^s_{p,q} (H)
\quad \text{and} \quad 
\big\| e^{-tH} f \big\|_{ B^{s}_{p,q}(H)} 
\leq C \| f \|_{ B^{s}_{p , q} (H) }
\end{equation} 
for all $t \geq 0$. 
\\
{\rm (ii)} If $s _2 \geq s_1$, $p_1 \leq p_2$ and 
$$
d \Big(\frac{1}{p_1} - \frac{1}{p_2} \Big) 
          + s_2 - s_1 > 0, 
$$
then there exists a constant $C > 0$ such that 
\begin{equation}\label{317-12}
\big\| e^{-tH} f \big\|_{ B^{s_2}_{p_2,q_2}(H)} 
\leq C t^{-\frac{d}{2} (\frac{1}{p_1} - \frac{1}{p_2}) 
          - \frac{s_2 - s_1}{2}
         } 
     \| f \|_{ B^{s_1}_{p_1 , q_1} (H) } 
\end{equation}
for any $f \in B^{s_1}_{p_1, q_1} (H)$. 
\end{thm}

We also have the continuity property of the semigroup in our Besov spaces 
associated with $H$ 
as well as the standard Besov spaces. 

\begin{thm}\label{thm:cts}
Let $s \in \mathbb R$, 
$1 \leq p,q \leq \infty$ and 
$1/p + 1/p' = 1/q + 1/q' = 1$. 
\\
{\rm (i)} Assume that $q < \infty$ and $f \in B^s_{p,q} (H)$. 
Then 
\begin{equation}\notag 
\lim _{t \to 0} \big\| e^{-t H} f -f 
\big\|_{ B^{s}_{p,q}(H)} = 0.
\end{equation}
{\rm (ii)} 
Assume that $1 < p \leq \infty$, $q = \infty$ and $f \in B^s_{p,\infty} (H)$. 
Then $e^{-t H} f$ converges to $f$ in the dual weak sense 
as $t \to 0$, namely, 
$$
\lim _{ t \to 0 } 
\sum _{j \in \mathbb Z} 
\int _{{\mathbb R^d}} 
   \Big\{ \phi_j (\sqrt{H})\big( e^{-tH} f -f \big)\, 
   \Big\}
   \overline{ \phi_j (\sqrt{H})g} \, dx 
= 0 
$$
for any $g \in B^{-s}_{p',1} (H)$. 

\end{thm}
\vskip3mm

We have an equivalent norm of the Besov spaces by using the semigroup. 

\begin{thm}\label{thm:equiv}
Let $s , s_0\in \mathbb R$, 
$s_0 > s/2$ and $1 \leq p,q \leq \infty$. 
Recall $j_0$, which was introduced in \eqref{0321-2} (i.e. $j_0=-1$). 
\black  
Then there exists a constant $C>0$ such that 
\begin{equation}\label{318-1}
C^{-1}\| f \|_{ B^s_{p,q} (H)} 
\leq 
\Big\{ \int _0^{2^{-2j_0}}
     \Big(  t^{-\frac{s}{2}} \| (tH ) ^{s_0} e^{-tH } f \|_{X} 
     \Big)^q
     \frac{dt}{t}
\Big\} ^{\frac{1}{q}}
\leq C \| f \|_{ B^s_{p,q} (H)} 
\end{equation}
for any $f \in  B^s_{p,q} (H)$, 
where $X$ can be $L^p$ or $B^0_{p,r} (H)$ with $1 \leq r \leq \infty$. 
\end{thm}

\noindent 
{\bf Remark. } Recalling $f_j = \phi_j(\sqrt{H})f$, 
we can estimate 
$\| (tH)^{s_0} e^{-tH} f_j \|_{L^p} $ 
by 
$ (t2^{2j})^{s_0} e^{-t2^{2j}} \| f_j \|_{L^p} $, 
which leads us to 
\[
\begin{split}
\Big\{ 
\int_0 ^\infty \Big( t^{-\frac{s}{2}} \| (tH)^{s_0} e^{-tH} f_j \|_{L^p} 
\Big) ^q \dfrac{dt}{t}\Big\}^{\frac{1}{q}}
\simeq 
&
\Big\{ 
\int_0 ^\infty \Big( t^{-\frac{s}{2}}
  (t2^{2j})^{s_0} e^{-t2^{2j}}
\Big) ^q \dfrac{dt}{t}\Big\}^{\frac{1}{q}}
\|  f_j \|_{L^p} 
\simeq 2^{sj} \| f_j \|_{L^p}. 
\end{split}
\]
On the other hand, the change of variable $t \mapsto 2^{-2j} t$ in the middle integral above implies 
\[
\begin{split}
2^{sj} \| f_j \|_{L^p} 
\simeq 
& 
\Big\{ 
\int_{2^{-2j}} ^{2^{-2(j-1)}} 
 \Big( t^{-\frac{s}{2}}
  (t2^{2j})^{s_0} e^{-t2^{2j}}
\Big) ^q \dfrac{dt}{t}\Big\}^{\frac{1}{q}}
\|  f_j \|_{L^p} . 
\end{split}
\]
Then, summing up in $j$ results in (\ref{318-1}). Since $j_0 =-1$ which is related to $\inf \sigma(H)>0$, the interval of the integral in the middle term of (\ref{318-1}) is only a bounded interval near $t=0$. 
We may see from this fact that our case corresponds to the inhomogeneous case of the Besov space for the standard Laplacian.

\black

\vskip3mm

The following theorem states the maximal regularity estimate for the semigroup. 

\vskip3mm

\begin{thm}\label{thm:max_reg}
Let $s \in \mathbb R$ and $1 \leq p,q \leq \infty$. 
Assume that $u_0 \in B^{s+2 - \frac{2}{q}}_{p,q} (H) $, 
$ f \in L^q (0,\infty ; B^s_{p,q}(H))$. Let $u $ be given by  
\begin{gather}
\notag 
u(t) = e^{- t H } u_0 
   + \int_0^t e^{-(t-\tau) H } f(\tau) d\tau . 
\end{gather}
Then there exists a constant $C>0$ independent of $u_0$ and $f$ such that 
\begin{equation}\label{318-2}
\| \partial _t u \|_{L^q (0,\infty ; B^s_{p,q}(H))}
+ \| H u \|_{L^q (0,\infty ; B^s_{p,q}(H))} 
\leq 
C\| u_0 \|_{ B^{s+2 - \frac{2}{q}}_{p,q} (H) }
+ C \| f \| _{L^q (0,\infty ; B^s_{p,q}(H))} . 
\end{equation}
\end{thm}

We finally mention a generalization of our results for the specific operator $H = - \Delta + |x|^2$ to more general Schr\"odinger operators with a potential $V$ that diverges at infinity, as studied in \cite{Y-2004}, where 
the potential $V$ is assumed to satisfy the following conditions for some $m > 2$ \black:
\begin{enumerate}
    \item[(a)] There exist constants $R > 0$ and $C_1 > 0$ such that 
    \[ \frac{1}{C_1} (1+|x|^2)^{\frac{m}{2}}\leq V(x) \leq C_1 (1+|x|^2)^{\frac{m}{2}} \quad \text{for } |x| \geq R.
    \]
    \item[(b)] For every multi-index $\alpha$, there exists a constant $C_\alpha > 0$ such that 
    \[ |\partial _x ^\alpha V(x)| \leq C_\alpha (1+|x|^2)^{\frac{m-|\alpha|}{2}} . 
    \]
\end{enumerate}
It is possible to introduce the Besov spaces associated with $-\Delta + V(x)$, as was done in \cite{IMT-2019}. We can then expect that the corresponding results stated in the introduction of this paper hold for these generalized operators. We also remark that the diverging property of the potential is crucial for showing the compact embedding in Proposition~\ref{prop:intro1} (v).

\black 

\black
\vspace{3mm}

This paper is organized as follows. Essentially, our tools for the bilinear estimates in Theorem \ref{thm:bilinear} rely on the Leibnitz rules applied to the operator $H$, and commutative properties with the multiplication by $x$ and the derivatives $\nabla$. We prepare some lemmas to describe such practical results in Section 2. Section 3 is devoted to the proof of Theorem \ref{thm:bilinear}. Since Theorems 1.2-1.5 
may be proved in a similar way in the existing literature, we will briefly add
explanations on the proofs in the Appendix.  \black

\vskip3mm 

\section{Preliminary}

In this section, we prepare some useful lemmas for the proof of Theorem 
\ref{thm:bilinear}.

\vskip3mm

\begin{prop} \label{prop:1} (\cite{dg}) 
For any $p\in (1,\infty)$ and $\alpha>0$, there exists a constant $C>0$ such that
\begin{equation*}
C^{-1} \|H^{\alpha} f\|_{L^p(\R^d)}
\leq \|(-\Delta)^{\alpha} f\|_{L^p(\R^d)} + \||x|^{2\alpha} f\|_{L^p} 
\le C \|H^{\alpha} f\|_{L^p(\R^d)}.
\end{equation*}
\end{prop}

\begin{lem}\label{lem:0805-2}For every multi-indices $\alpha , \beta$, 
there exists a constant $C>0$ such that 
\begin{equation}\label{0803-1}
\| x^\alpha \nabla ^\beta f  \|_{L^2}
\leq C \| H^{\frac{|\alpha|+|\beta|}{2}}  f  \|_{L^2}
\end{equation}
for 
all $f\in L^2$ satisfying $H^{\frac{|\alpha|+|\beta|}{2}}  f \in L^2$. 
\black
\end{lem}
\begin{proof}
If $\alpha = (0,\cdots,0)$ or $\beta = (0,\cdots, 0)$, 
then Proposition~\ref{prop:1} proves the inequality \eqref{0803-1}. 
It is sufficient to prove the case when $\alpha \not= 0$ and  $\beta \not = 0$. 
Also it is sufficient to prove for $f \in \mathcal S(\mathbb R^d)$ by the density argument. 
\black

When $|\alpha|= |\beta | = 1$, we estimate  
$x_j \partial _{x_k} f$ ($j,k = 1,2, \cdots , d$), 
\begin{eqnarray*}
\| x_j \partial _{x_k} f \|_{L^2}^2
&=& \int x_j^2 |\partial _{x_k} f|^2 dx 
\leq \int (H \partial _{x_k}f ) \overline{\partial_{x_k} f} dx .
\end{eqnarray*}
Since $H\partial_{x_k} = \partial_{x_k} H -2 x_k$, 
we have 
\begin{eqnarray*}
\| x_j \partial _{x_k} f \|_{L^2}^2
&\leq &  
\int ( \partial _{x_k} Hf ) \, \overline{\partial_{x_k} f} dx 
-2 \int x_k f \, \overline{\partial_{x_k}f}dx
\\
& \leq & \| Hf \| _{L^2}\| \Delta f \|_{L^2} + 2 \| x_k f \|_{L^2} \| \nabla f \|_{L^2}.
\end{eqnarray*}
It follows by Proposition~\ref{prop:1} and 
$\inf \sigma (H) > 0$ 
that 
$$\| x_k f \|_{L^2}, \| \nabla f \|_{L^2} 
\leq C \| H^{\frac12} f \|_{L^2} \leq C \| Hf \|_{L^2},$$
thus we obtain 
\[
\| x_j \partial _{x_k} f \|_{L^2}^2
\leq C \| H f \|_{L^2}^2.
\]

We apply the induction argument for the proof of the higher order cases. 
Let $M \geq 2$ be a natural number and we assume that 
\[
\| x^\alpha \nabla ^\beta f\|_{L^2} 
\leq C \| H^{\frac{|\alpha|+|\beta|}{2}} f \|_{L^2} \quad 
\text{ if } |\alpha| + |\beta| \leq M.
\]
Let us prove the estimate when $|\alpha| + |\beta| = M + 1$. 

If $|\alpha|$ is an even number, then by Proposition \ref{prop:1} \black
\[
\| x^\alpha  \nabla ^\beta f\|_{L^2} 
\leq C \| H^\frac{|\alpha|}{2} \nabla ^\beta f \|_{L^2}. 
\]
with $\frac{|\alpha|}{2} \in \N$. \black
Since $H \nabla = \nabla H -2x$,  
there exist $\Lambda_M$ a subset consisting of indices $(\alpha' , \beta')$ 
with the total order $|\alpha'|+|\beta'|$ less than $M$
and positive constants $C_{\alpha' , \beta'}$ such that 
\\ 
\[
H^\frac{|\alpha|}{2} \nabla ^\beta f
= \nabla ^\beta H^\frac{|\alpha|}{2} f 
+ \sum_{(\alpha' , \beta ') \in \Lambda_{M}} 
C_{\alpha',\beta'} \,  x^{\alpha'} \nabla ^{\beta'} f,
\]
which proves that 
\[
\| x^\alpha  \nabla ^\beta f\|_{L^2} 
\leq C \| \nabla ^\beta H^\frac{|\alpha|}{2} f \|_{L^2}
+ C \sum_{(\alpha', \beta') \in \Lambda _M}
    \| x^{\alpha'}  \nabla ^{\beta'} f\|_{L^2} .  
\]
Proposition~\ref{prop:1} and the assumption of the induction imply that 
\[
\| x^\alpha  \nabla ^\beta f\|_{L^2} 
\leq C \| H^\frac{|\alpha|+|\beta|}{2} f \|_{L^2}
+ C \sum_{(\alpha', \beta') \in \Lambda _M}
    \| H^\frac{|\alpha'|+|\beta'|}{2} f \|_{L^2} .
\]
We also know  $\inf \sigma(H) > 0$ \black and obtain 
the inequality \eqref{0803-1}. 

If $|\alpha|$ is an odd number and $|\alpha| \geq 3$, then we write 
\[
x^\alpha = x_j x^{\tilde \alpha} , \quad \nabla ^\beta = \partial_{x_k} \nabla ^{\tilde \beta}, 
\quad \text{for some } \tilde \alpha, \tilde \beta,  j,k
\quad \text{where $|\tilde \alpha|$ is even}, 
\]
and by the integration by parts, 
\begin{eqnarray*}
\| x^\alpha  \nabla ^\beta f\|_{L^2} ^2 
&=& \int x_j^2 x^{\tilde \alpha} \nabla ^\beta f \cdot 
   \overline{x^{\tilde \alpha}\nabla ^\beta f} dx
\\
&=& -\int x_j^2 x^{\tilde \alpha} \nabla ^{\tilde \beta} f \cdot 
   \overline{x^{\tilde \alpha} \partial_{x_k} \nabla ^\beta f} dx
- \int \Big( \partial_{x_k} (x_j^2 x^{\tilde \alpha} x^{\tilde \alpha}) \Big) \nabla ^{\tilde \beta} f \cdot 
   \overline{ \nabla ^\beta f} dx .
\end{eqnarray*}

\noindent 
\black
Since $2 + |\tilde \alpha|$ and $|\tilde \alpha|$ are even, 
we have by the Cauchy Schwarz inequality and the previous argument 
for even number polynomials that 
\[
\bigg| 
\int x_j^2 x^{\tilde \alpha} \nabla ^{\tilde \beta} f \cdot 
   \overline{x^{\tilde \alpha} \partial_{x_k} \nabla ^\beta f} dx
\bigg|
\leq C \| H^{1+\frac{|\tilde \alpha|}{2}} \nabla^{\tilde\beta} f \|_{L^2}
\| H^{\frac{|\tilde \alpha|}{2}} \partial_{x_k}\nabla^{\beta} f \|_{L^2}
\leq C \| H^{\frac{|\alpha|+|\beta|}{2}} f \|_{L^2}^2. 
\]
For the second term, we notice 
$|\tilde \beta| = |\beta|-1$ and that 
the order of the polynomial 
$\partial _{x_k} (x_j^2 x^{\tilde \alpha } x^{\tilde \alpha})$ 
is $1 + 2 |\tilde \alpha| = 2 |\alpha| -1$ at most, and 
we apply the multiplication of the polynomials of order $|\alpha|$ and $|\alpha|-1$ 
by $\nabla ^{\tilde \beta} f $ and $ \overline{\nabla ^{\beta} f}$, 
respectively. 
We then write 
\[
\Big( \partial _{x_k} (x_j^2 x^{\tilde \alpha } x^{\tilde \alpha})\Big) 
\nabla ^{\tilde \beta} f \cdot \overline{\nabla ^{\beta} f}
= \sum_{(\alpha ' , \beta') \in \Lambda '_{M-1}} 
  x^{\alpha'} \nabla ^{ \beta'} f 
  \sum_{(\alpha '' , \beta'') \in \Lambda ''_{M-1}} 
    \overline{x^{\alpha''}\nabla ^{\beta''} f},
\]
where $\Lambda '_{M-1}, \Lambda ''_{M-1}$ are sets of multi-indices 
for polynomials and derivatives such that the sum of the two orders 
are $M-1$. 
We apply the assumption of the induction for $M-1$ to have that 
\[
\bigg| \int \Big( \partial_{x_k} (x_j^2 x^{\tilde \alpha} x^{\tilde \alpha}) \Big) \nabla ^{\tilde \beta} f \cdot 
   \overline{ \nabla ^\beta f} dx
\bigg|
\leq C \| H^{\frac{M-1}{2}} f \|_{L^2}^2 . 
\]
The above two inequalities proves the case when $|\alpha| + |\beta| = M+1$, 
and we conclude the estimate \eqref{0803-1}. 
\black 
\end{proof}

The following lemma is fundamental for our argument
and will be used several times in our proof. 
It is the uniform boundedness of the spectral multiplier
with derivatives and multiplication by polynomials.

\begin{lem}\label{lem:0805-1} 
For multi-indices $\alpha$ and $\beta$, there exists a positive constant $C_{\alpha, \beta}$ \black
such that for every $ f \in L^p $ and $j \in \mathbb Z$, 
\begin{equation}
\| x^\alpha \nabla ^\beta \phi_j(\sqrt{H}) f  \|_{L^p}
+ \| \phi_j(\sqrt{H}) x^\alpha \nabla ^\beta f  \|_{L^p}
\leq C_{\alpha , \beta } 2^{ (|\alpha|+|\beta|)j} 
  \| f \|_{L^p}. 
\end{equation}
\black 
\end{lem}

Let us give a comment on the proof of this lemma. 
In the case where $p=1$, we can apply the lemmas below
and the argument in \cite{IMT-RMI} to the operator with derivatives and polynomials 
to prove the inequality. 
The case where $p=\infty$ follows from the duality argument, and the case where $1 < p < \infty$ is proved by interpolation.

\vskip3mm

To prove Lemma~\ref{lem:0805-1}, 
we introduce a set $\mathscr A_N$ of some bounded operators on $L^2 (\mathbb R^d)$ 
and scaled amalgam spaces $\ell ^1 (L^2)_{\theta}$ for $\theta > 0$ 
to prepare a lemma. 
Hereafter, for $k \in \mathbb Z^d$, 
$C_\theta (k)$ denotes a cube with the center $\theta^{\frac{1}{2}} k$ 
and side length $\theta ^{\frac{1}{2}}$, namely, 
$$
C_{\theta} (k) 
:= \big\{ x \in \mathbb R^d \, \big| \,
    |x_j - \theta^{\frac{1}{2}} k_j| 
    \leq 2^{-1} \theta^{\frac{1}{2}} \text{ for } j = 1,2, \cdots, d
    \, \big\}, 
$$ 
and $\{ \chi _{C_\theta (k)} \}_{k \in \mathbb Z^d} 
\subset C_0^\infty (\mathbb R^d)$ is a partition of the 
unity such that 
\[
\chi_{C_1 (k)} (x) = \chi_{C_1(0)}(x-k), \quad 
\chi_{C_\theta (k)}(x) = \chi_{C_1(k)}(\theta ^{-\frac{1}{2}}x),
\]
\[
\sum_{k \in \mathbb Z^d} \chi_{C_\theta (k)} (x) = 1 \text{ for all } x \in \mathbb R^d. 
\]

\noindent 
{\bf Definition. }  
For $N \in \mathbb N$, 
$\mathscr{A}_N$ denotes the set of all 
bounded operators $T$ on $L^2 (\Omega)$ such that 
$$
\| T \|_{\mathscr{A}_N} := 
\sup_{k \in \mathbb{Z}^d} 
 \big\| |\cdot - \theta^{1/2}k|^{N} T \chi_{C_{\theta}(k)}
 \big\|_{L^2 \to L^2} < \infty. 
 $$

\vskip3mm 

\noindent {\bf Remark.} 
We remark that this partition of the unity consists of smooth functions, 
while in 
the reference~\cites{IMT-RMI} the authors use non-smooth functions to compose a partition of the unity. 
We need some smoothness of the partition 
to study the operators $T = x^\alpha \nabla ^\beta \varphi(\sqrt{H})$ 
and $\varphi(\sqrt{H}) x^\alpha \nabla ^\beta$. 
 The lemmas (Lemma \ref{lem:317-2} and Lemma \ref{lem:317-0}) below hold 
also for our partition of the unity \black and is proved with suitable modification, 
but we omit the detail. 

\vskip3mm 

\noindent 
{\bf Definition. } 
The space $\ell ^1 (L^2)_{\theta}$ is defined by letting  
$$
\ell ^1 (L^2)_{\theta} := 
\big\{ f \in L^2 _{\rm loc} (\mathbb R^d) \, \big| \, 
   \| f \|_{\ell^1 (L^2)_{\theta}} < \infty
  \, \big\}, 
$$
where
$$
\displaystyle 
\| f \|_{\ell^1 (L^2)_{\theta}} 
:= \sum _{k \in \mathbb Z^d} 
   \| f \|_{L^2 (C_\theta (k))}. 
$$

\begin{lem} \label{lem:317-0} {\rm (\cite{IMT-2019})}
The operator $\varphi(\sqrt{\theta H})$ with $\varphi \in \mathcal{S}(\R)$ belongs to $\mathscr{A}_N$ for any $\theta >0$. Moreover,  
there exists a constant $C>0$ such that 
$$
\|\varphi(\sqrt{\theta H})\|_{\mathscr{A}_N} \le C \theta^{\frac{N}{2}}.
$$
\end{lem}
\black
\begin{lem}\label{lem:317-2} {\rm (\cite{IMT-RMI})}
\noindent 
{\rm (i)} 
Let $N \in \mathbb N$ and $N > d/2$. Then 
there exists a constant $C > 0$ such that 
\begin{equation}
\label{317-6}
\| T \|_{\ell^1 (L^2)_{\theta} \to \ell^1 (L^2)_{\theta}} 
\leq C 
  \Big( \| T \|_{L^2 \to L^2} 
    + \theta ^{-\frac{d}{4}} 
      \| T \|_{\mathscr A _N} ^\frac{d}{2N} 
      \| T \|_{L^2 \to L^2} ^{ 1- \frac{d}{2N}}
  \Big) 
\end{equation}
for any $T \in \mathscr A_N$ and $\theta>0 $. 

\noindent 
{\rm (ii)} 
Let $\beta$ be a real number satisfying $ \beta > d /4$.  Then there exists a constant $C > 0$ such that 
\begin{equation}
\label{317-8}
\big\| (1 + \theta H) ^{-\beta} \big\| _{L^1 \to \ell^1 (L^2)_{\theta}} 
\leq C \theta^{ -\frac{d}{4} } 
\end{equation}
for any $\theta > 0$. 
\end{lem}

\begin{lem}
Let $m, N \in \mathbb N $. For every multi-indices $\alpha, \beta$ 
with $|\alpha|+|\beta| \leq 2m$, 
there exists a constant $C > 0$ such that 
\begin{equation}
\label{317-7-1}
\|x^\alpha \nabla ^\beta \psi ( (1 + \theta H)^{-m} ) \|_{\mathscr A_N} 
\leq C \theta ^{\frac{N}{2}-\frac{|\alpha|+|\beta|}{2}} , 
\end{equation}
\begin{equation}
\label{317-7-2}
\|\psi ( (1 + \theta H)^{-m} )x^\alpha \nabla ^\beta  \|_{\mathscr A_N} 
\leq C \theta ^{\frac{N}{2}-\frac{|\alpha|+|\beta|}{2}} . 
\end{equation}
for any $\psi \in \mathcal S (\mathbb R)$ 
with ${\rm supp \,} \psi \subset [0,\infty)$ 
and $\theta > 0$. 
\end{lem}

\begin{proof}
We write $R _\theta = (1+ \theta H)^{-1}$ and notice that 
\[
\begin{split}
(x- \theta^{\frac{1}{2}} n) x^\alpha \nabla ^\beta 
=& x^\alpha \nabla ^\beta (x-\theta ^\frac{1}{2}n)
 + x^\alpha (x \nabla ^\beta - \nabla ^\beta x) 
\\
=& x^\alpha \nabla ^\beta (x-\theta ^\frac{1}{2}n)
 + 
 \begin{cases}
 x^\alpha \nabla ^{\tilde \beta} , \text{ with } 
 |\tilde \beta| = |\beta| -1,
 \\
 \text{or } \quad 0 ,
 \end{cases}
\end{split}
\]
by which we can write for $j = 1, 2, \cdots, d$ 
\[
\begin{split}
& (x_j- \theta^{\frac{1}{2}}n_j)^N x^\alpha \nabla ^\beta \psi(R_\theta^m) 
\\
=& x^{\alpha} \nabla ^\beta (x_j-\theta^{\frac{1}{2}}n_j)^N \psi(R_\theta^m) 
 + x^\alpha \sum_{( \beta ', k) \in \Lambda _{|\beta|-1, N-1}} 
   C_{\beta', k} \nabla ^{\beta'} (x_j - \theta^{\frac{1}{2}}n_j)^k \psi (R_\theta^m),
\end{split}
\]
where $\Lambda _{|\beta|-1, N-1}$ is a subset of 
multi-indeces $(\beta',k)$ such that $|\beta'|\leq |\beta|-1$ 
and $k \leq N-1$, 
and we need to study the $L^2$-boudedness. 
Here we focus on the first term above, 
since the total orders of the derivatives and the polynomials are less 
and the second term can be handled similarly to the first term. 
We then consider the commutator for the first term, 
\[
\begin{split}
&x^{\alpha} \nabla ^\beta (x_j-\theta^{\frac{1}{2}}n_j)^N \psi(R_\theta^m)
\\
=& 
x^{\alpha} \nabla ^\beta \psi(R_\theta^m)(x_j-\theta^{\frac{1}{2}}n_j)^N 
+ x^{\alpha} \nabla ^\beta 
  \Big( (x_j-\theta^{\frac{1}{2}}n_j)^N \psi(R_\theta^m)
     - \psi(R_\theta^m)(x_j-\theta^{\frac{1}{2}}n_j)^N 
  \Big) 
\\
=& I + II, 
\end{split}
\]
and the first term is handled by Lemma~\ref{lem:0805-2}, 
\begin{equation}\label{0806-1}
\begin{split}
& 
\displaystyle  \sup_{n \in \mathbb Z^d}  
 \|I \cdot \chi_{C_\theta (n)}\|_{L^2 \to L^2} 
\\
 \leq
 & \| x^\alpha \nabla^\beta H^{-\frac{|\alpha|+|\beta|}{2}} \|_{L^2\to L^2}
 \| H^{\frac{|\alpha|+|\beta|}{2}} \psi(R_\theta ^m) \|_{L^2 \to L^2}
 \| (x_j-\theta^{\frac{1}{2}}n_j)^N \cdot \chi_{C_\theta (n)} \|_{L^2\to L^2}
\\
\leq & C \theta ^{-\frac{|\alpha|+|\beta|}{2}} \theta ^{\frac{N}{2}} 
\sup_{\lambda > 0} \lambda^{\frac{|\alpha|+|\beta|}{2}} \psi ((1+\lambda)^{-m}),
\end{split}
\end{equation}
where the supremum above is finite because of 
$\psi(0)= 0$ and smooth at the origin. 
Next, we recall the formula
\[
\psi(R_\theta^m) = (2\pi)^{-\frac{1}{2}} 
 \int_{-\infty}^\infty e^{it R_\theta^m} \widehat{\psi} (\xi) dt ,
\]
from which it is sufficient to study 
\begin{eqnarray*}
&&\displaystyle  \sup_{n \in \mathbb Z^d}  
 \Big\|
 \int_{-\infty}^\infty 
  x^\alpha \nabla ^\beta 
  \Big(  \Big( (x_j-\theta^{\frac{1}{2}}n_j)^N e^{itR_\theta^m}
     - e^{itR_\theta^m} (x_j-\theta^{\frac{1}{2}}n_j)^N 
  \Big) 
  \Big) 
 \widehat{\psi}(t) dt \cdot  \chi_{C_\theta (n)}
 \Big\|_{L^2 \to L^2} .  
\end{eqnarray*}
We follow the argument with commutators (see the proof of Lemma~6.3 in~\cite{IMT-RMI}), but we only explain the different point. 
The problem for the commutator is reduced to $R_\theta^m$ instead of 
$R_\theta$ in~\cite{IMT-RMI}, since 
\[
\begin{split}
&(x_j-\theta^{\frac{1}{2}}n_j) e^{itR_\theta^m}
     - e^{itR_\theta^m} (x_j-\theta^{\frac{1}{2}}n_j)
\\
=& 
\int_0^t \partial _s \Big( 
 e^{i(t-s)R_\theta^m} (x_j-\theta^{\frac{1}{2}}n_j) e^{isR_\theta^m}
\Big)  ds 
\\
=& 
i\int_0^t 
 e^{i(t-s)R_\theta^m} 
 \Big( -R_\theta^m (x_j-\theta^{\frac{1}{2}}n_j) 
      +  (x_j-\theta^{\frac{1}{2}}n_j) R_\theta ^m 
 \Big) e^{isR_\theta^m}
\Big)  ds 
\\
=& 
i\int_0^t 
 e^{i(t-s)R_\theta^m} R_\theta^m 
 \Big( -x_j (1+\theta H)^m + (1+\theta H)^m x_j
 \Big) R_\theta ^m  e^{isR_\theta^m}
\Big)  ds . 
\end{split}
\]
We then see that 
there exist $\Lambda_{2k-1}$ a subset of $\alpha', \beta'$ with the total 
order $|\alpha'| + |\beta'| \leq 2k-1$ ($k = 1, \cdots, m$) 
and constants $C_{\alpha', \beta'}$ such that 
\[
-x_j (1+\theta H)^m + (1+\theta H)^m x_j
= \sum_{k=1}^{m} \theta ^k \sum_{(\alpha',\beta')\in \Lambda_{2k-1}} 
C_{\alpha', \beta'} x^{\alpha '} \nabla ^{\beta'} . 
\]
We can then handle by the boundedness of 
$x^\alpha \nabla ^\beta R_\theta ^m, 
x^{\alpha'} \nabla ^{\beta '} R_\theta ^m$ in $L^2$ 
proved by Lemma~\ref{lem:0805-1}. 
Therefore we conclude that 
\[
\|x^\alpha \nabla ^\beta \psi ( R_\theta ^m ) \|_{\mathscr A_N} 
\leq C\theta ^{-\frac{|\alpha|+|\beta|}{2}} \theta^{\frac{N}{2}} . 
\]

We explain how to prove the second inequality \eqref{317-7-2} by a similar argument 
to the proof above.  We write 
\[
\begin{split}
& (x_j-\theta^{\frac{1}{2}}n_j)^N \psi(R_\theta^m)x^{\alpha} \nabla ^\beta
\\
=& 
\psi(R_\theta^m) x^{\alpha} \nabla ^\beta (x_j-\theta^{\frac{1}{2}}n_j)^N
+ \Big( (x_j-\theta^{\frac{1}{2}}n_j)^N \psi(R_\theta^m) x^{\alpha} \nabla ^\beta 
     - \psi(R_\theta^m)x^{\alpha} \nabla ^\beta (x_j-\theta^{\frac{1}{2}}n_j)^N 
  \Big), 
\end{split}
\]
and apply the boundedness of the operators 
$R_\theta ^m x^\alpha \nabla ^\beta $ in $L^2$, which is proved 
by Lemma~\ref{lem:0805-2} and the duality argument provided that 
$|\alpha| + |\beta| \leq 2m$. 
In fact, 
we may have 
\[
\big\| 
(x_j-\theta^{\frac{1}{2}}n_j)^N \psi(R_\theta^m)x^{\alpha} \nabla ^\beta (\chi_{C_\theta (n)} f)
  \big)
\big\|_{L^2}
\leq C\theta ^{-\frac{|\alpha|+|\beta|}{2}} \theta^{\frac{N}{2}} 
 \| f \|_{L^2}, 
\]
for $f \in \mathcal S (\mathbb R^d)$, 
where the above constant $C$ is independent of 
$n$ and $f$. 
A density argument implies that 
\[
\big\| 
(x_j-\theta^{\frac{1}{2}}n_j)^N \psi(R_\theta^m)x^{\alpha} \nabla ^\beta \chi_{C_\theta (n)} 
  \big)
\big\|_{L^2 \to L^2}
\leq C\theta ^{-\frac{|\alpha|+|\beta|}{2}} \theta^{\frac{N}{2}} 
\]
for all $n \in \mathbb Z^d$ and the constant $C$ is 
independent of $n$. 
We then obtain the second inequality \eqref{317-7-2}. 
\end{proof}

\begin{proof}[Proof of Lemma~\ref{lem:0805-1}]
As explained below Lemma~\ref{lem:0805-1}, we only 
prove the case when $p=1$. 
We also introduce $\theta$ such that $\theta = 2^{-2j}$. 

We write by the partition of the unity $\{\chi_{C_\theta(n)} \}_{n \in \mathbb Z^d}$,  
\[
\| x^\alpha \nabla ^\beta \phi _j(\sqrt{H}) f \|_{L^1} 
\leq  \sum _{n \in \mathbb Z^d} 
\|  \chi_{C_\theta(n)} x^\alpha \nabla ^\beta \phi _j(\sqrt{H}) f \|_{L^1} 
\leq C \theta ^{\frac{d}{4}} 
  \| x^\alpha \nabla ^\beta \phi _j(\sqrt{H}) f \|_{l^1L^2(C_\theta(n))} . 
\]
Given a positive real number $\gamma>\frac{d}{4}$,  \black we choose 
$\widetilde \varphi \in \mathcal S(\mathbb R)$ as  
\[
\widetilde \varphi(\lambda) = (\lambda + 1)^{  \black} \phi_j \black (\sqrt{\lambda}), 
\quad \lambda > 0 . 
\]
By the definition of $\widetilde \varphi$, we have 
\[
\theta ^{\frac{d}{4}} 
  \| x^\alpha \nabla ^\beta \phi _j(\sqrt{H}) f \|_{l^1L^2(C_\theta(n))} 
= \theta ^{\frac{d}{4}} 
  \| x^\alpha \nabla ^\beta \widetilde \varphi(\theta H) (\theta H +1)^{-\gamma} \black f \|_{l^1L^2(C_\theta(n))} . 
\]
It follows from \eqref{317-6}, \eqref{0803-1} and 
\eqref{317-7-1} that for $T = x^\alpha \nabla ^\beta \widetilde \varphi(\theta H) $
\begin{eqnarray*}
&& 
\theta ^{\frac{d}{4}} 
  \| x^\alpha \nabla ^\beta \widetilde \varphi(\theta H) (\theta H +1)^{- \gamma \black} f \|_{l^1L^2(C_\theta(n))} 
\\
&\lesssim & \theta ^{\frac{d}{4}} 
 \Big( \| T \|_{L^2 \to L^2} 
    + \theta ^{-\frac{d}{4}} 
      \| T \|_{\mathscr A _N} ^\frac{d}{2N} 
      \| T \|_{L^2 \to L^2} ^{ 1- \frac{d}{2N}}
  \Big) 
  \|(\theta H +1)^{- \gamma \black} f \|_{l^1L^2(C_\theta(n))} 
\\
&\lesssim&\theta ^{\frac{d}{4}} 
 \Big( \theta ^{-\frac{|\alpha|+|\beta|}{2}}
    + \theta ^{-\frac{d}{4}} 
      (\theta ^{\frac{N}{2}-\frac{|\alpha|+|\beta|}{2}}) ^\frac{d}{2N} 
      (\theta ^{-\frac{|\alpha|+|\beta|}{2}}) ^{ 1- \frac{d}{2N}}
  \Big) 
  \|(\theta H +1)^{- \gamma \black} f \|_{l^1L^2(C_\theta(n))} 
\\
&\lesssim&
  \theta ^{\frac{d}{4}}  \cdot 
  \theta ^{-\frac{|\alpha|+|\beta|}{2}}
  \|(\theta H +1)^{- \gamma \black} f \|_{l^1L^2(C_\theta(n))} . 
\end{eqnarray*}
We finally apply \eqref{317-8} and have that 
for $\theta = 2^{-2j}$
\[
  \theta ^{\frac{d}{4}}  \cdot 
  \theta ^{-\frac{|\alpha|+|\beta|}{2}}
  \|(\theta H +1)^{-\gamma \black} f \|_{l^1L^2(C_\theta(n))} . 
\lesssim 
  \theta ^{\frac{d}{4}}  \cdot 
  \theta ^{-\frac{|\alpha|+|\beta|}{2}}
  \cdot \theta ^{-\frac{d}{4}} 
  \| f \|_{L^1} 
= 2^{(|\alpha|+|\beta|)j} \| f \|_{L^1} ,
\]
which proves Lemma~\ref{lem:0805-1}. 
\end{proof}

\black

\section{Proof of Theorem~\ref{thm:bilinear}}

In this section, we give a proof for Theorem ~\ref{thm:bilinear}. 
Note that the derivatives and multiplications of functions are taken
in the $\mathcal{S}'$ sense. \black

We start with the proof of Proposition 1.2, item (i).  
For each $j \in \mathbb Z$, we write 
\[
\begin{split}
\phi_j(\sqrt{H}) (f\ol g) 
=& 
 \phi_j(\sqrt{H})
 \Big( \sum_{l: |l-j|\leq 1} + \sum _{l: l-j \leq -2} + \sum_{l: l-j \geq 2} 
 \Big)
\Big( \sum_{k: k \leq l-2}f_k g_l \Big).
\end{split}
\]
We can handle the first case $|l-j| \leq 1$ in the same way as in standard
Besov spaces associated with the Laplacian (see, for example, \cite{bcd}). However, for the sake of completeness, we explain briefly here. 
For $j\in \mathbb Z$ given, using the boundedness of spectral multiplier \black 
 (\ref{unifbdd})  
and the H\"{o}lder inequality, we get 
\begin{eqnarray*}
2^{sj}\| \phi_j(\sqrt{H}) \sum_{l:|l-j|\le 1} \sum_{k:k\le l-2} f_k g_l\|_{L^p} 
&\le & 2^{sj} 
\sum_{l:|l-j|\le 1} \|g_l\|_{L^{p_2}} 
\|\sum_{k:k\le l-2} f_k \|_{L^{p_1}}\\
& \lesssim &  \black  \sum_{l:|l-j|\le 1} 2^{sl} \|g_l\|_{L^{p_2}} 
\|f\|_{L^{p_1}} 
\end{eqnarray*}
since
\begin{equation} \label{sumbdd}
\Big \|\sum_{k:k\le l-2} f_k \Big \|_{L^{p_1}} \le \|f\|_{L^{p_1}}.
\end{equation}
In fact, by introducing $\psi _{l-2} = \sum_{k \leq l-2} \phi_k$, 
we can apply the uniform bound \eqref{unifbdd} to $\psi_{l-2}$ instead of $\phi_j$ to obtain (\ref{sumbdd}). \black
We then take the $l^q$ norm and apply the Young inequality. 
$$
\left[\sum_{j\in \mathbb Z} 
\left\{ 2^{sj}\| \phi_j(\sqrt{H}) \sum_{l:|l-j|\le 1} \sum_{k:k\le l-2} f_k g_l\|_{L^p}\right \}^q  \right]^{\frac{1}{q}} 
\le C\|f\|_{L^{p_1}} \|g\|_{B^{s}_{p_2,q}} \left(\sum_{|j|\le 1} 2^{-sj}\right) $$
where the sum $\sum_{|j|\le 1} 2^{-sj}$ is finite. We point out that in the case of standard Laplacian $-\Delta$ we do not need to consider the case $|l-j|\ge 2$ since the supports of decomposition functions are disjoint, but in our case we need it.   
In this proof below, we will see that even if we consider $H=-\Delta+|x|^2$, since this case can be treated as a {\it perturbation} from the Laplacian $-\Delta$ case, the same bilinear estimates follow-namely the term of the case $|l-j|\ge 2$ should be small. 
Such an approach is inspired \black by the argument presented in \cite{dp}, where the equivalence between the two Besov spaces with and without a potential is discussed.

Let us consider the second case 
$l-j \leq -2$.  Take $m$ with $2m > |s|$ and fix. We see that  
\[
\phi_j(\sqrt{H}) \Big( \sum_{k \leq l-2}f_k g_l \Big)
= \phi_j(\sqrt{H}) H^{-m} H^m\Big( \sum_{k \leq l-2}f_k g_l \Big) ,
\]
and it follows by the uniform boundedness of the spectral multiplier 
\eqref{unifbdd} and the Leibniz rule with Lemma~\ref{lem:0805-1} (we do not use the equivalence of the Sobolev norm here because it excludes the cases $p=1,\infty$) that 
for each $l$
\begin{eqnarray} \label{0409-1-1}
\Big\| \phi_j(\sqrt{H}) H^{-m} H^m \black  \Big( \sum_{k \leq l-2}f_k g_l \Big) \Big\|_{L^p} 
\leq
C 2^{-2mj} 
   \Big\| H^m\Big( \sum_{k \leq l-2}f_k g_l \Big) \Big\|_{L^p}.
\end{eqnarray}
Here we remark that 
\begin{eqnarray}\label{0409-1-2}
 2^{-2mj} 
   \Big\| \sum_{k \leq l-2}f_k H^m g_l \Big\|_{L^p}
   \leq 
C_m 2^{-2m j}\| f \|_{L^{p_1}} 2^{2m l}\| g_l \|_{L^{p_2}}. 
\end{eqnarray}
In this inequality (\ref{0409-1-2}), we 
focused on the most important term $f_k H^m g_l$ in the right hand side of (\ref{0409-1-1}), and 
\black 
applied \eqref{sumbdd} for the term 
$\sum_{k \leq l-2} f_k$.\footnote{Writing $g_l = (\phi_{l-1}(\sqrt{H}) + \phi_{l}(\sqrt{H}) + \phi_{l+1}(\sqrt{H})) g_l$, we apply the spectral multiplier theorem to 
$(\phi_{l-1}(\sqrt{H}) + \phi_{l}(\sqrt{H}) + \phi_{l+1}(\sqrt{H}))$, 
and we can then keep $g_l$ in the inequality.\black} 

For the other terms, $f_k$ should have multiplication by the polynomials 
or derivatives of order one at least, and it follows by $2^k \leq 2^l$ that
\begin{equation}\notag 
\begin{split}
 2^{-2mj} 
\Big\| H^m\Big( \sum_{k \leq l-2}f_k g_l \Big) 
  - \sum_{k \leq l-2}f_k  H^m g_l\Big\|_{L^p} 
  \leq 
  & C 
 2^{-2mj}  
 \sum _{k \leq l-2} 2^k \| f_k \|_{L^{p_1}} 2^{(2m-1)l} \| g_l \|_{L^{p_2}}
\\
  \leq
  &C 2^{-2mj} 
   \| f \|_{L^{p_1}} 2^{2ml} \| g_l \|_{L^{p_2}}.
\end{split}
\end{equation}
We apply this inequality, replace $l$ by $l + j$, and then  
the second case $l-j \le -2$ \black can be estimated by 
\[
\begin{split}
& \Big\{ \sum_{j \in \mathbb Z}  
  \Big( 2^{sj} 2^{-2mj} \sum_{l-j \leq -2} \| f \|_{L^{p_1}} 2^{2ml} \| g_l \|_{L^{p_2}}
  \Big)^q
  \Big\}^{\frac{1}{q}}
\\
\leq& C\Big\{ \sum_{j \in \mathbb Z}  
  \Big( 2^{sj} 2^{-2mj} \sum_{l \leq -2} \| f \|_{L^{p_1}} 2^{2m(l+j)} \| g_{l+j} \|_{L^{p_2}}
  \Big)^q
  \Big\}^{\frac{1}{q}}
\\
\leq& C
\sum_{l \leq -2} 2^{(2m-s)l}\| f \|_{L^{p_1}}
\Big\{ \sum_{j \in \mathbb Z}  
  \Big( 2^{s(l+j)} \| g_{l+j} \|_{L^{p_2}}
  \Big)^q
  \Big\}^{\frac{1}{q}}
\\
\leq & C \| f \|_{L^{p_1}} \| g \|_{B^s_{p_2,q}(H)}. 
\end{split}
\]
Finally we treat the third case $l-j \geq 2$. We observe that the quantity  $\phi_j(\sqrt{H}) \Big( \sum_{k \leq l-2}f_k g_l \Big)$ can be rewritten as follows.  
\\

\begin{eqnarray} \nonumber
&&\phi_j(\sqrt{H}) \Big( \sum_{k \leq l-2}f_k g_l \Big)
= \phi_j(\sqrt{H}) \Big( \sum_{k \leq l-2}f_k \cdot  H^m H^{-m} g_l \Big) 
\\ \label{eq:leibnitz2}
&=& \sum_{(\alpha' , \alpha'', \beta', \beta'') \in \Lambda_{2m}} 
C_{\alpha' , \alpha'', \beta', \beta''}
\phi_j(\sqrt{H})x^{\alpha'}\nabla ^{\beta'} 
\Big\{ \Big( \sum_{k \leq l-2} x^{\alpha''}\nabla ^{\beta''}f_k \Big) \cdot  H^{-m}g_l \Big\},
\end{eqnarray}
where $\Lambda_{2m}$ is a subset of indices such that 
$|\alpha'|+ |\alpha''| + |\beta'| + |\beta''| \leq 2m$; 
indeed, when $m=1$, we can write as follows.
\begin{eqnarray*}
\phi_j(\sqrt{H}) \Big( \sum_{k \leq l-2}f_k H H^{-1} g_l \Big)
&=& 
\phi_j(\sqrt{H}) H
 \Big( \sum_{k \leq l-2}f_k H^{-1}g_l \Big) 
+ 2\phi_j(\sqrt{H}) 
 \Big( \sum_{k \leq l-2} (-\Delta  f_k ) H^{-1} g_l \Big)
\\
&& 
+2  \phi_j(\sqrt{H}) \nabla \cdot 
 \Big( \sum_{k \leq l-2} (\nabla f_k)  H^{-1} g_l \Big). 
\end{eqnarray*}
Then we assume that (\ref{eq:leibnitz2}) holds for $m$, and prove that the case $m+1$ holds. Indeed, we write the $m+1$ case as 
\begin{eqnarray*}
\phi_j(\sqrt{H}) \Big( \sum_{k \leq l-2}f_k \cdot  H^{m+1} H^{-(m+1)} g_l \Big) 
=\phi_j(\sqrt{H}) \Big( \sum_{k \leq l-2}f_k \cdot  H H^{-1} [H^m H^{-m} g_l] \Big)
\end{eqnarray*}
and use the result for $m=1$ case. Then we have
\begin{eqnarray*}
\phi_j(\sqrt{H}) \Big( \sum_{k \leq l-2}f_k H H^{-1} [H^m H^{-m} g_l] \Big)
&=& 
\phi_j(\sqrt{H}) H
 \Big( \sum_{k \leq l-2}f_k H^{-1}[ H^m H^{-m} g_l] \Big) \\
&& + 2\phi_j(\sqrt{H}) 
 \Big( \sum_{k \leq l-2} (-\Delta  f_k ) H^{-1} [H^m H^{-m}g_l] \Big)
\\
&& 
+2  \phi_j(\sqrt{H}) \nabla \cdot 
 \Big( \sum_{k \leq l-2} (\nabla f_k)  H^{-1} [H^m H^{-m} g_l] \Big). 
\end{eqnarray*}
Now we apply (\ref{eq:leibnitz2}) with $m$ for $f_k \mapsto f_k, -\Delta f_k, \nabla f_k$ respectively and   
$g_l \mapsto H^{-1}g_l$. We may then see that the case $m+1$ also holds for 
(\ref{eq:leibnitz2}). 
\black

\vskip3mm

Now, in the third case $l-j \ge 2$, \black dividing the sum in $k \le l-2$ \black into two cases $k < j$, $k \geq j$,  
Lemma~\ref{lem:0805-1} yields that 

\begin{eqnarray*}
\Big\| \phi_j(\sqrt{H}) \Big( \sum_{k \leq l-2}f_k g_l \Big)\Big\|_{L^p} 
&\leq&  
\Big \|\phi_j(\sqrt{H})
\sum_{k<j, k\le l-2} f_k g_l \Big\|_{L^p}
+
\Big\| \phi_j (\sqrt{H}) \sum_{j \leq k \leq l-2} f_k g_l \Big\|_{L^p} \\
&\le &
C  
2^{2mj}\| f \|_{L^{p_1}} 2^{-2ml} \|g_l\|_{L^{p_2}} 
+ 
C
\sum_{j \leq k \leq l-2} 2^{2mk} \| f_k \|_{L^{p_1}} 2^{-2ml}\| g_l \|_{L^{p_2}}.
\end{eqnarray*}
Indeed, the first term has been estimated as follows. 
\begin{eqnarray*}
&&\Big\| \phi_j(\sqrt{H}) \Big( \sum_{ k<j, k \leq l-2}f_k g_l \Big)\Big\|_{L^p} 
\\
&\leq&  
\Big\|\sum_{(\alpha' , \alpha'', \beta', \beta'') \in \Lambda_{2m}} 
C_{\alpha' , \alpha'', \beta', \beta''}
\phi_j(\sqrt{H})x^{\alpha'}\nabla ^{\beta'} 
\Big\{ \Big( \sum_{k <j, k \leq l-2} x^{\alpha''}\nabla ^{\beta''}f_k \Big) \cdot  H^{-m}g_l \Big\} \Big\|_{L^p} \\
&\leq& 
\sum_{(\alpha' , \alpha'', \beta', \beta'') \in \Lambda_{2m}} 
C_{\alpha' , \alpha'', \beta', \beta''}
2^{(|\alpha'|+|\beta'|)j}
\Big\| \sum_{k <j, k \leq l-2} x^{\alpha''}\nabla ^{\beta''}f_k \Big\|_{L^{p_1}} \|  H^{-m}g_l \|_{L^{p_2}}\\
&\leq& 
\sum_{(\alpha' , \alpha'', \beta', \beta'') \in \Lambda_{2m}} 
C_{\alpha' , \alpha'', \beta', \beta''}
2^{(|\alpha'|+|\beta'|)j} 
2^{(|\alpha''|+|\beta''|)k}
\Big\|\sum_{k <j, k \leq l-2} f_k \Big\|_{L^{p_1}} \|  H^{-m}g_l \|_{L^{p_2}}\\
&\leq& C2^{2mj} \|f \|_{L^{p_1}} \|  H^{-m}g_l \|_{L^{p_2}}, 
\end{eqnarray*}
where we have used Lemma \ref{lem:0805-1} \black  twice.

\black

This implies that for the case of $k < j$,
\[
\begin{split}
& \Big\{ \sum_{j \in \mathbb Z}  
  \Big( 2^{sj}  \sum_{l-j \geq -2} 2^{2mj }  \| f \|_{L^{p_1}} 
    2^{-2ml} \| g_l \|_{L^{p_2}}
  \Big)^q
  \Big\}^{\frac{1}{q}}
\\
\leq & \| f \|_{L^{p_1}}  
  \Big\{ \sum_{j \in \mathbb Z}  
  \Big( 2^{(s+2m)j}  \sum_{l \geq -2}  
    2^{-2m(l+j)} \| g_{l+j} \|_{L^{p_2}}
  \Big)^q
  \Big\}^{\frac{1}{q}}
\\
\leq & \| f \|_{L^{p_1}}  \sum_{l \geq -2}  2^{-(s+2m)l}
  \Big\{ \sum_{j \in \mathbb Z}  
  \Big( 
    2^{s(l+j)} \| g_{l+j} \|_{L^{p_2}}
  \Big)^q
  \Big\}^{\frac{1}{q}}
\\
\leq & C \| f \|_{L^{p_1}} \| g \|_{ B^s_{p_2,q}(H)},
\end{split}
\]
and for the case of $k \geq j$, 
\[
\begin{split}
& \Big\{ \sum_{j \in \mathbb Z}  
  \Big( 2^{sj}  \sum_{l-j \geq -2} \sum_{j \leq k \leq l-2} 2^{2mk }  \| f_k \|_{L^{p_1}} 
    2^{-2ml} \| g_l \|_{L^{p_2}}
  \Big)^q
  \Big\}^{\frac{1}{q}}
\\
= & \Big\{ \sum_{j \in \mathbb Z}  
  \Big( 2^{sj}  \sum_{l \geq -2} \sum_{0 \leq k \leq l-2} 2^{2m(k+j) }  \| f_{k+j} \|_{L^{p_1}} 
    2^{-2m(l+j)} \| g_{l+j} \|_{L^{p_2}}
  \Big)^q
  \Big\}^{\frac{1}{q}}
\\
\leq & C \| f \|_{L^{p_1}}
\sum_{l \geq -2} 2^{(-s-2m)l}\sum_{0 \leq k \leq l-2} 2^{2mk }
\Big\{ \sum_{j \in \mathbb Z}  
  \Big( 2^{s(l+j)}    \| g_{l+j} \|_{L^{p_2}}
  \Big)^q
  \Big\}^{\frac{1}{q}}
\\
\leq & C \| f \|_{L^{p_1}}
\Big( \sum_{l \geq -2} 2^{-sl} \Big) 
\| g \|_{ B^s_{p_2,q}(H)} . 
\end{split}
\]
\vspace{3mm}

The $l\le k-2$ case, i.e. the product rule for $f \og g $ can be shown in the same way. Further, 
the item (ii) is also similarly proved. 
\vspace{3mm}



\black

Next we show (iii). First, as above 
we decompose for $j \ge -1$ fixed, 
$$ \phi_j(\sqrt{H}) (f\oo g) = 
\phi_j(\sqrt{H}) \left\{\sum_{k:k\ge j-2} \sum_{|l-k|\le 1} f_k g_l
+ \sum_{k:k\le j-2} \sum_{l:|l-k|\le 1}f_k g_l \right\} =(\mathrm{I})+(\mathrm{II})$$
We first consider (I). By (\ref{unifbdd}), 
\begin{eqnarray*}
2^{(s_1+s_2)j} \|\phi_j(\sqrt{H}) (\mathrm{I}) \|_{L^p} 
&=&2^{(s_1+s_2)j} \Big\|\phi_j(\sqrt{H}) \sum_{k:k\ge j-2} \sum_{|l-k|\le 1} f_k g_l \Big\|_{L^p} \\
&\le&  2^{(s_1+s_2)j} \sum_{k:k\ge j-2} \sum_{|l-k|\le 1}
\|f_k \|_{L^{p_1}} \|g_l\|_{L^{p_2}} \\
&\le&  \sum_{k:k\ge j-2} 
2^{-(k-j)(s_1+s_2)} 2^{ks_1} \|f_k\|_{L^{p_1}} 
\sum_{|\nu|\le 1} 2^{(k-\nu) s_2} 2^{s_2 \nu}\|g_{k-\nu}\|_{L^{p_2}}.
\end{eqnarray*}
Thus we take the $l^q$ norm in $j$, and use the Young inequality and 
the H\"{o}lder inequality with $\frac{1}{q}=\frac{1}{q_1}+\frac{1}{q_2}$ to conclude 
$$
\Big\{\sum_j \{ 2^{(s_1+s_2)j} \|\phi_j(\sqrt{H}) (\mathrm{I}) \|_{L^p} \}^q\Big\}^{\frac{1}{q}} 
\le  \Big(\sum_{j \ge -2} 2^{-(s_1+s_2)j} \Big) \|f\|_{B^{s_1}_{p_1,q_1}}
\|g\|_{B^{s_2}_{p_2,q_2}},
$$
where $\sum_{j \ge -2} 2^{-(s_1+s_2)j} <\infty$ if $s_1+s_2>0$.
Next, for the term (II), 
as above, take 
$m$ such that $2m>s_1+s_2$. 
For a fixed $j \ge -1$, 
we write 


\[
\phi_j (\sqrt{H}) \sum_{k:k\le j-2} 
\sum_{l:|l-k|\le 1} f_k g_l
= H^{-m} \phi_j (\sqrt{H}) \sum_{k:k\le j-2} 
\sum_{l:|l-k|\le 1} H^{m} (f_k g_l),
\]
and apply a similar argument as in \eqref{0409-1-1} 
with the condition that $|l-k| \leq 1$. 
By the Leibniz rule, we write 
\[
H^m (f_k g_l) 
= \sum _{(\alpha , \beta , \alpha' , \beta ') \in \Lambda _{2m}}
 C_{\alpha , \beta, \alpha ', \beta '} 
 (x^\alpha \nabla ^\beta f_k )
 (x^{\alpha '} \nabla ^{\beta '} g_l ) , 
\]
where $\Lambda _{2m} $ 
is a set of multi-indices for polynomials and derivatives such that 
the total order $|\alpha| + |\beta| + |\alpha'| + |\beta'|$ 
is less than $2m$. 
We apply Lemma~\ref{lem:0805-1} and see that 
\[
\| (x^\alpha \nabla ^\beta f_k )
 (x^{\alpha '} \nabla ^{\beta '} g_l )\|_{L^p} 
 \leq C 
 2^{(|\alpha| + |\beta|)k} \| f_k \|_{L^{p_1}}
 2^{(|\alpha'| + |\beta|')l} \| g_l \|_{L^{p_2}},
\]
where the multi-indices must satisfy
\[
|\alpha| + |\beta| \leq 2m , \quad 
|\alpha'| + |\beta'| \leq 2m, \quad 
|\alpha| + |\beta| + |\alpha'| + |\beta'| \leq 2m .
\]
We can then write 
\begin{eqnarray*}
&&2^{(s_1+s_2)j} \Big\|\phi_j (\sqrt{H}) \sum_{k:k\le j-2} 
\sum_{l:|l-k|\le 1} f_k g_l \Big\|_{L^p} 
\\ 
&\lesssim &
2^{(s_1+s_2)j} \sum_{k:k\le j-2} 2^{-2mj} 
\sum_{l:|l-k|\le 1} \sum _{\substack{0\le n,\tilde{n} \le 2m, \\ 0\le n+\tilde{n} \le 2m}}
 C_{n,\tilde{n}} 2^{n k} 
\|f_k \|_{L^{p_1}} 
\cdot 
2^{\tilde{n}l} \|g_{k-\nu}\|_{L^{p_2}} 
\\
&\le &
2^{(s_1+s_2)j} \sum_{k:k\le j-2} 2^{-2mj} 
\sum_{l:|l-k|\le 1} \sum_{n=0}^{2m} 
\sum_{\tilde{n}=0}^{2m-n} 
 C_{n,\tilde{n}} 2^{n k} 
\|f_k \|_{L^{p_1}} 
\cdot 
2^{\tilde{n}l} \|g_{k-\nu}\|_{L^{p_2}} 
\\
&\le &
2^{(s_1+s_2)j} \sum_{k:k\le j-2} 2^{-2mj} 
\sum_{l:|l-k|\le 1} \sum_{n=0}^{2m} 
 C'_{n,m} 2^{n k} 
\|f_k \|_{L^{p_1}} 
\cdot 
2^{(2m-n)l} \|g_{k-\nu}\|_{L^{p_2}} 
\\
&\leq &  2^{(s_1+s_2)j} \sum_{k:k\le j-2} 2^{-2mj} 
\sum_{|\nu|\le 1} \sum _{n=0}^{2m}
 C'_{m,n} 2^{n k} 
\|f_k \|_{L^{p_1}} 
\cdot 
2^{(2m-n) (k-\nu)} \|g_{k-\nu}\|_{L^{p_2}} 
\\ 
&=& 
2^{(s_1+s_2)j} \sum_{k:k\le j-2} 2^{-2mj} 
\sum_{|\nu|\le 1} \sum_{ n=0}^{2m}
2^{(2m - n) \cdot (-\nu)} 
C'_{m, n} 2^{2mk} 
\|f_k \|_{L^{p_1}} \|g_{k-\nu}\|_{L^{p_2}} \black  \\
&\le & \sum_{k:k\le j-2} 2^{(2m-(s_1+s_2))(k-j)} 2^{s_1 k} \|f_k\|_{L^{p_1}} \\
&& \times \sum_{|\nu|\le 1} 2^{s_2 \nu} \sum_{n=0}^{m} C'_{m,n} 2^{s_2 (k-n\nu)} \|g_{k-\nu}\|_{L^{p_2}}  ,
\end{eqnarray*}
where $C'_{m,n}$ ($n = 0, 1,2, \cdots, 2m$) are appropriate real numbers 
depending on the subscript $m, n$. 
\\
\black

\noindent 
Now, first, we take $l^q$ norm in $j$, then use the Young inequality, 
\begin{eqnarray*}
\Big\| \sum_{k:k\le j-2} \sum_{l:|l-k|\le 1} f_k g_l \Big\|_{B^{s_1+s_2}_{p,q}} &\le& 
\sum_{j\le -2} 2^{(2m-(s_1+s_2))j} \|f\|_{B^{s_1}_{p_1,q_1}} \\
&& \times \sum_{n=0}^{m} C_{m,n} \Big\|\sum_{|\nu|\le 1} 2^{{s_2} \nu} 
2^{s_2 (\cdot-\nu)} 2^{s_2 (1-n) \nu } \|g_{\cdot-\nu}\|_{L^{p_2}}  \Big\|_{l^{q_2}}
\end{eqnarray*}
Again using the Young inequality in the last term we get
$$\Big\| \sum_{k:k\le j-2} \sum_{l:|l-k|\le 1} f_k g_l \Big\|_{B^{s_1+s_2}_{p,q}}
\le 
C_{|s_2|,m} \|f\|_{B^{s_1}_{p_1,q_1}}\|g\|_{B^{s_2}_{p_2,q_2}}. 
$$
Combining (i)-(iii), we obtain (iv). 
\hfill\qed

\vskip3mm

\appendix

\section{proof of (\ref{unifbdd})}

In this section, we give a brief proof for the uniform bound 
(\ref{unifbdd}).

\begin{proof}[Proof of (\ref{unifbdd})] 
It is sufficient to show $L^1$-estimate for 
$\varphi( \sqrt{\theta H})$, then $L^\infty$-estimate follows by the 
duality argument. We then can make use of 
the Riesz-Thorin interpolation theorem 
to obtain $L^p$-estimates for $1\le p\le \infty$. 

Recalling the definition 
$l^1(L^2)_{\theta}$ in Section 2, 
we obtain 
\begin{align*}
\|\varphi(\sqrt{ \theta  H}) f\|_{L^1}
& = 
\sum_{n\in\Z^d}
\|
\varphi(\sqrt{\theta  H}) f\|_{L^1(C_{\theta}(n))}\\
& \le
\sum_{n\in\Z^d}
|C_\theta(n)|^{1/2} 
\|
\varphi( \sqrt{\theta  H}) f\|_{L^2(C_{\theta}(n))}\nonumber\\
& \le 
\theta^{d/4} 
\|
\varphi( \sqrt{ \theta H}) f\|_{l^1(L^2)_{\theta}}, \nonumber
\end{align*}
where we have used the bound
$
|C_\theta(n)|^{1/2}\le \theta^{d/4}.
$
\vspace{3mm}

For $\beta>0$, we consider 
$\tilde{\varphi} \in \mathcal S (\mathbb R)$ defined by
\[
\tilde{\varphi} (\lambda)
=
(\lambda + M)^\beta \varphi (\lambda) \quad 
\text{for }
\lambda \in \sigma(H), \quad M>-1. 
\]

Note that we may write 
\begin{align*}
\|\varphi( \sqrt{\theta H}) f\|_{l^1(L^2)_{\theta}}
=&
\left\|
\tilde{\varphi}(\sqrt{\theta H}) 
(\theta H + M )^{-\beta} f 
\right\|_{l^1(L^2)_{\theta}}. 
\end{align*}
Using Lemma \ref{lem:317-2}, we get 
\begin{align*}
&\left\|
\tilde{\varphi}( \sqrt{\theta H})
(\theta H + M )^{-\beta} f 
\right\|_{l^1(L^2)_{\theta}}\\
\le & 
C \Big( 
\|\tilde{\varphi}( \sqrt{\theta H})\|_{L^2 \to L^2} 
+ \theta^{-d/4} 
{\|{\tilde{\varphi}
(\sqrt{\theta H})}}\|^{d/2N}_{\mathscr{A}_N} 
\|\tilde{\varphi}( \sqrt{\theta H})
\|^{1-d/2N}_{L^2\to L^2} 
\Big)\\
& \times 
\left\|
(\theta H + M )^{-\beta}  f\right\|_{l^1(L^2)_{\theta}}.
\end{align*}

Remark that the bound in $L^2$ follows from 
$$ 
\|\tilde{\varphi}( \sqrt{\theta H})f\|_{L^2} 
\le \int_{\inf \sigma(H)}^{\infty} |\tilde{\varphi} (\sqrt{\theta \lambda})|^2 
d \|E(\lambda) f\|_{L^2}^2 \le \|\tilde{\varphi}\|_{L^{\infty}}\|f\|_{L^2}^2
$$
for any $ \theta > 0.$
Moreover, thanks to Lemma \ref{lem:317-0}, 
the right hand side is estimated by 
\begin{align*}
C \left\{ 1 + \theta^{-d/4} \cdot 
(\theta^{N/2})^{d/2N} \right\} \theta^{-d/4} 
\|f\|_{L^1}
=
C \theta^{-d/4} \|f\|_{L^1}, 
\end{align*}
provided $\beta$ 
satisfies $\beta >d/4$. 
Summarizing those estimates, 
we find that 
\begin{equation*}
\|\varphi(\sqrt{\theta  H}) f\|_{l^1(L^2)_{\theta}}
\le 
C \theta^{-d/4} \|f\|_{L^1}.
\end{equation*}
Therefore, we conclude that
\[
\|\varphi (\sqrt{\theta H})f\|_{L^1}
\le
C\|f\|_{L^1}
\]
for any $\theta > 0 $ and $f\in L^1$.
\end{proof}

\black
\section{proof of Proposition \ref{prop:intro1} (v) (vi) and (vii)}
\label{appendix:B}

\begin{proof}[Proof of Proposition \ref{prop:intro1} (v)] 
We begin by proving the continuous embedding $L^p \hookrightarrow B^s_{p,1} (H)$ for $s < 0$. This embedding is a fundamental result in the theory of non-homogeneous Besov spaces, and the proof proceeds as follows.
Since $s < 0$, we have by (\ref{unifbdd}) \black
\[
\| f \|_{B^s_{p,1}(H)} 
\leq \sum _{j =-1}^\infty 2^{sj} \| \phi_j(\sqrt{H})  f\|_{L^p}
\leq \Big(\sum _{j =-1}^\infty 2^{sj} \Big)  \cdot C \| f \|_{L^p}. 
\]
\black

We prove the compact embedding of 
$B^\alpha_{p,\infty}(H) \subset L^p$ 
provided that $\alpha>0$. 

Let $\{ f_n \}_{n=1}^\infty$ be a bounded sequence in 
$B^\alpha _{p,\infty}$. 
For $R > 0$, we have 
\[
\sup _{ n \in \mathbb N} \| f_n \|_{L^p (\{ |x| >R \})}
\leq R^{-\frac{\alpha}{2}}
\sup _{ n \in \mathbb N} \| |\cdot|^\frac{\alpha}{2} f_n \|_{L^p (\{ |x| >R \})}. 
\]
\[
\| |\cdot|^\frac{\alpha}{2} f_n \|_{L^p (\{ |x| >R \})} 
\leq 
\sum _{j = -1}^\infty  \| |\cdot|^\frac{\alpha}{2} \phi_j(\sqrt{H}) f_n \black \|_{L^p}
\leq C \sum _{j = -1}^\infty 2^{\frac{\alpha}{2}j -\alpha j} 
 \cdot 2^{\alpha j} \| \phi_j(\sqrt{H}) f_n \black \|_{L^p},
\]
which implies that 
\begin{equation}\notag 
\sup _{ n \in \mathbb N} \| f_n \|_{L^p (\{ |x| >R \})}
\leq C R^{-\frac{\alpha}{2}}
\sup _{ n \in \mathbb N} \| f_n \|_{B^\alpha_{p,\infty}(H)} . 
\end{equation}

Define 
\[
f_{n,J} =\sum_{j=-1}^J \phi_j(\sqrt{H})f_n.
\]
Since $f_{n,J}$ is smooth, we have from the Arzel\'a-Ascoli theorem that 
for each $J$, $M$, a subsequence $\{ f_{n_k,J, M} \}_{k=1}^\infty $ exists 
such that it converges uniformly to an continuous function, 
which we denote by $f^M_{\leq J}$, 
on each compact set $\{|x| \le M\}$ of $\mathbb R^d$. 
We may choose a subsequence satisfying the monotonicity 
$\{ f_{n_k,J, M+1} \}_{k=1}^\infty 
\subset \{ f_{n_k,J, M} \}_{k=1}^\infty$ with respect to $M$, 
and we then find  
$\{ f_{n_k, J} \}_{k=1}^\infty$ such that 
\[
\{ f_{n_k, J} \}_{k=1}^\infty 
\subset \{ f_{n_k,J, M+1 \black} \}_{k=1}^\infty  
\quad \text{for all } M ,
\]
and $\{ f_{n_k, J} \}_{k=1}^\infty $ converges 
uniformly to the function $f^M_{\leq J}$ on 
the compact set $\{ |x| \leq M \}$ for each $M$ and $J$. 
By the monotonicity with respect to $M$ and the uniform convergence, 
\[
f_{\leq J}^M = f_{\leq J }^{M+1} \quad \text{in }  \{ |x| \leq M \} , 
\]
for all $M$. We introduce the function $f_{\leq J}$ such that 
\[
f_{\leq J} (x) = f_{\leq J}^M (x) , \quad \text{ if } |x| \leq M , 
\quad M = 1,2, \cdots. 
\]
and  according to above arguments we see that $f_{n_k, J}$ converges to 
$f _{\leq J}$ uniformly on $\{ |x| \leq M \}$ 
for each $M \in \mathbb N$. 

\black 
We easily see from Fatou's lemma for $p < \infty$ 
and an elementary argument for $p = \infty$ that $f_{\leq J} \in L^p$. 
Moreover, 
\[
\| \phi_j (\sqrt{H})f_{\leq J} \|_{L^p} 
\leq \liminf_{k \to \infty} \| \phi_j(\sqrt{H}) f_{n_k,J}  \|_{L^p}.
\]
and we may also have 
\begin{equation}\notag 
\begin{split}
\phi_j(\sqrt{H}) f_{\leq J} 
=&  \phi_j (\sqrt{H}) f_{\leq J+J'}
\quad \text{ for all } J' = 1,2, \dots.  
\end{split}
\end{equation}
Let us define 
\[
f = \sum _{j=-1}^\infty \phi_j(\sqrt{H}) f_{\leq j+1}.
\]
 We then notice that 
\begin{equation}\notag 
\begin{split}
&\| f \|_{B^\alpha_{p,\infty}(H)} 
\leq 
 \liminf_{k \to \infty} \|  f_{n_k,J}  \|_{B^\alpha _{p,\infty}(H)}, 
 \quad  \text{for each } J
\\
&  \| |\cdot|^{\frac{\alpha}{2}} f \|_{L^p (\{ |x| >R \})}
\leq C R^{-\frac{\alpha}{2}} \| f \|_{B^\alpha_{p,\infty}(H)} 
\leq C R^{-\frac{\alpha}{2}} \sup _{n \in \mathbb N}\| f_n \|_{B^\alpha_{p,\infty}(H)} . 
\end{split}
\end{equation}

We take an arbitrary positive number $\varepsilon$. 
For each $J$, we consider 
the subsequence $\{ n_k \}_{k=1}^\infty$  associated with $f_{n_k,J}$, 
and the norm of $f_{n_k} -f$. 
We write  
\[
\begin{split}
\| f_{n_k} -f \|_{L^p} 
\leq & 
\| f_{n_k} -f_{n_k,J} \|_{L^p} + \| f_{n_k,J} -f \|_{L^p} 
\\
\leq 
& \| f_{n_k} -f_{n_k,J} \|_{L^p} 
+ \| f_{n_k,J} - f_{\le J \black} \|_{L^p (\{ |x| \leq M \})}
\\
& 
+ \| f_{\le J \black} -  f\|_{L^p(\{|x| \leq  M \})}
+ \| f_{n_k,J} -f \|_{L^p(|x| > M)} .
\\
\end{split}
\]
The first term is the norm for the high-spectrum part, 
and we have 
\[
\| f_{n_k} - f_{n_k,J} \|_{L^p} 
\leq \sum _{j \geq J-1} 2^{-\alpha j} \cdot 
2^{\alpha j} \| \phi_j (\sqrt{H})f_{n_k} \|_{L^p}
\leq C 2^{-\alpha J} 
 \sup _{n \in \mathbb N} \| f_{n}\|_{B^\alpha_{p,\infty}}. 
\]
The third term is the norm for the high-spectrum part, 
and we have 
\[
\| f_{\le J \black} - f \|_{L^p(|x|\le M)\black} 
\leq \sum _{j \geq J-1} 2^{-\alpha j} \cdot 
2^{\alpha j} \| \phi_j (\sqrt{H})f \|_{L^p}
\leq C 2^{-\alpha J} 
 \| f \|_{B^\alpha_{p,\infty}}. 
\]
The fourth term is bounded by
\[
\| f_{n_k,  J\black} - f\|_{L^p(\{ |x| > M \})} 
\leq C M ^{-\frac{\alpha}{2}} 
\sup_{n \in \mathbb N} \| f_{n}\|_{B^\alpha_{p,\infty}} . 
\]

We here choose $J$ such that the first and the third terms 
are small, i.e., 
\[
\|f_{n_k} - f_{n_k,J}  \|_{L^p} + \| f_{\le J} - f \|_{L^p(|x|\le M)} 
\leq 
 C 2^{-\alpha J} 
 \sup_n \| f_n \|_{B^\alpha_{p,\infty}} \black 
< \varepsilon, 
\]
and we may find $M$ such that
the fourth term can be small as follows. 
\[
\| f_{n_k, J \black} - f\|_{L^p(\{ |x| > M \})} 
\leq 
2C M ^{-\frac{\alpha}{2}} 
\sup_{n \in \mathbb N} \| f_n \|_{B^\alpha_{p,\infty}}
< 2\varepsilon. 
\]
We know that $f_{n_k,J}$ converges to $f_{\leq J}$ uniformly on 
$\{ |x| \leq M \}$ as $k \to \infty$, 
and then see that 
a natural number $ k_0(J,M)$ exists such that for $k \geq k_0(J,M)$
\[
\| f_{n_k ,J} - f_{\le J \black}\|_{L^p \{ |x| \leq M \}} < \varepsilon ,
\]
which implies that 
\[ \| f_{n_k} -f  \|_{L^p} < 4\varepsilon . 
\]
Therefore, we can find a subsequence of $\{ f_{n} \}_{n=1}^\infty$ 
such that it converges to $f$ in $L^p$. 

\black 

We turn to prove the compact embedding from 
$B^\alpha _{p,\infty}(H)$ to $B^s_{p,1}(H)$ provided that $s < \alpha$. 
By the lifting property, we can assume that 
$s < 0 < \alpha$. 
Let $\{ f_n \}_{n=1}^\infty $ be a bounded sequence 
in $B^\alpha_{p,\infty} (H)$. 
By the previous proof, we can find a subsequence of $\{ f_n \}_{n=1}^\infty$ 
which converges in $L^p$, and it is easy to see that the subsequence 
converges in $B^s_{p,1} (H)$ by the continuous embedding 
$L^p \hookrightarrow B^s_{p,1} (H)$ provided $s < 0$, 
as we showed at the beginning of the proof. 
\end{proof}

\begin{proof}[Proof of Proposition \ref{prop:intro1} (vi)]
Since $\sum _{j \in \mathbb Z} \phi_j(\sqrt{H})$ 
is the identity operator, we have 
\begin{eqnarray*}
\|f\|_{L^p} &=& 
 \sum_{j\in \mathbb Z} \|\phi_j(\sqrt{H}) f\|_{L^p}
= \|f\|_{B^0_{p,1}}.
\end{eqnarray*}
On the other hand, 
$$\|f\|_{B^0_{p,\infty}}=\sup_{j\ge -1} \|\phi_j(\sqrt{H})f\|_{L^p} 
\le C \|f\|_{L^p},$$
by \eqref{unifbdd}. 
\end{proof}

\begin{proof}[Proof of Proposition \ref{prop:intro1} (vii)]
We follow the argument in Section 4 in \cite{FI}. 
To explain the idea, we only consider the case when $r ,r_0 \leq p $. 
We notice that 
\begin{equation}\notag 
s- \frac{d}{p} \in \Big( - \frac{d}{r} ,~ s_0 - \frac{d}{r_0} \Big) .
\end{equation}

We here recall the following inequality for $p \geq r$. 
\[
\| \phi_j(\sqrt{H}) \|_{L^p \to L^r} 
\leq C 2^{d(\frac{1}{r}-\frac{1}{p})j}, \quad \text{ for all } j \in \mathbb Z,
\]
which is a generalization of the boundedness \eqref{unifbdd} 
and we refer to Theorem~1.1 in \cite{IMT-RMI} 
(see also the proof of Theorem~1.2 (ii) below). 

For $N \in \mathbb Z$ we split the infinite series in the definition of the norm of Besov spaces into two series. 

\[
\begin{split} 
\| f \|_{B^s_{p,1}(H)} 
 \leq 
& C \sum _{ j \leq N} 2^{sj + d (\frac{1}{r}- \frac{1}{p})j}  
   \| \phi_j(\sqrt{H})f \|_{L^{r}} 
+ C\sum _{ j > N} 2^{sj + d(\frac{1}{r_0}-\frac{1}{p})j} 
  \| \phi_j(\sqrt{H}) f \|_{L^{r_0}} 
\\
\leq 
& C 2^{sN + d (\frac{1}{r}- \frac{1}{p})N}  
   \| f \|_{B^0_{r,\infty}(H)} 
+  C 2^{sN + d(\frac{1}{r_0}-\frac{1}{p})N - s_0N} 
\| f \|_{B^{s_0}_{r_0,\infty}(H)},
\end{split}
\]
since $s+d(1/r - 1/p) > 0, s + d(1/r_0 - 1/p) -s_0 < 0$. Choosing $N$ such that 
\[
 2^{sN + d (\frac{1}{r}- \frac{1}{p})N}  
   \| f \|_{B^0_{r,\infty}(H)} 
\simeq  2^{sN + d(\frac{1}{r_0}-\frac{1}{p})N - s_0N} 
\| f \|_{B^{s_0}_{r_0,\infty}(H)},
\]
we obtain the inequality in (vii). 
\black 
\end{proof}

\noindent 
{\bf Remark. } 
It is possible to prove a simpler inequality.  
\[
\| f \|_{B^\alpha _{p,q}} 
\leq \| f \|_{B^{\alpha _0}_{p_0 , q_0}} ^{1-\theta} 
  \| f \|_{B^{\alpha_1}_{p_1,q_1}} ^\theta
\]
where $\alpha , \alpha_0, \alpha _1 \in \mathbb R$, 
$\theta \in (0,1)$, $1 \leq p,p_0,p_1, q, q_0, q_1 \leq \infty$,
\[
\alpha = (1-\theta) \alpha_0 + \theta \alpha _1 , \quad 
\dfrac{1}{p} = \dfrac{1-\theta}{p_0} + \dfrac{\theta}{p_1}, \quad 
\dfrac{1}{q} = \dfrac{1-\theta}{q_0} + \dfrac{1}{q_1}. 
\]
In fact,  by the H\"{o}lder inequality, we get
$$\|f\|_{L^p} =\|f^{1-\theta+\theta}\|_{L^p} \le \|f\|^{1-\theta}_{L^{p_0}} \|f\|^{\theta}_{L^{p_1}}.$$
Therefore, 
$$2^{\alpha k}\|\phi_k(\sqrt{H})f\|_{L^p} 
\le 2^{\alpha_0 k (1-\theta)}\|\phi_k(\sqrt{H})f\|^{1-\theta}_{L^{p_0}} 
\hspace{1mm} 2^{\alpha_1 k\theta}\|\phi_k(\sqrt{H})f\|^{\theta}_{L^{p_1}}.$$
Then we take $l^q$ norm, and apply the H\"{o}lder inequality. 
\black

\section{proofs of Theorems 1.2-1.5}

Since the proof of Theorems~1.2--1.5 \black 
follows from the argument in \cite{Iw-2018}, we highlight only a few key points.
To this end, we prepare a lemma, which is similar to Lemma~2.1 in \cite{Iw-2018} 
for the Dirichlet Laplacian. 
Instead of the Dirichlet Laplacian, 
we consider the Hermite operator in this paper. 

\begin{lem}\label{lem:0321-1}
Let $N > d/2$, $1 \leq p \leq \infty$, $\delta > 0$ and $a,b > 0$. 
Then there exists a positive constant $C$, 
which depends on $N, \delta, a, b$, such that 
for any $\phi \in C_0 ^\infty (\mathbb R)$ with  
${\rm supp \, } \phi \subset [a,b] $, 
$G \in C^\infty((0,\infty ))$
\begin{equation}\label{317-1}
\| G(\sqrt{H}) \phi ( 2^{-j} \sqrt{H})  f \|_{L^p} 
\leq C \| G( 2^j \sqrt{\cdot } \, ) \phi (\sqrt{\cdot} \, ) 
       \|_{H^{N+\frac{1}{2} + \delta} (\mathbb R)} \| f \|_{L^p}
\end{equation}
for all $j \in \mathbb Z$, 
where 
$\| \phi \|_{H^s (\mathbb R) }
= \| (1+|\xi|^2)^{\frac{s}{2}} \mathcal F[\phi] \|_{L^2 (\mathbb R)}$. 
\end{lem}
\begin{proof}[Proof of Lemma~\ref{lem:0321-1}] 
The proof is similar to that of Lemma~2.1 in \cite{Iw-2018}, 
as the semigroup generated by the Hermite operator satisfies resolvent estimates
and the following Gaussian upper bound 
(see, e.g., Proposition~3.1 in \cite{IMT-RMI}). 
There exists a positive constant $C$ such that 
\begin{equation}\label{0326-1}
0 \leq e^{-t H} (x,y) 
\leq C t^{-\frac{d}{2}}e^{-\frac{|x-y|^2}{Ct}}, 
\text{ for all } t > 0 \text{ and } x,y \in \mathbb R^d,
\end{equation}
where $e^{-t H} (x,y)$ denotes the kernel of the operator $e^{-tH}$. 
We also refer to Section 6 in \cite{IMT-RMI} for results on Schr\"odinger operators, including the Hermite operator.
\end{proof}

\begin{proof}[Proof of Theorem 1.2 \black]
(i) It is well known that the kernel of the semigroup $e^{-tH}$ 
satisfies the Gaussian upper bound \eqref{0326-1}. 
This implies $L^p$ boundedness for all $1 \leq p \leq \infty$, 
which in turn proves boundedness of $e^{-tH}f$ \black in Besov spaces.
\\
(ii) We consider only the case $q_2 = 1$ and $q_1 = \infty$, 
as the embedding properties of Besov spaces allow us to deduce the other cases from this.

We define $\Phi_j = \phi_{j-1} + \phi_j + \phi_{j+1}$ and write
\[
\phi_j(\sqrt{H}) e^{-tH} f 
= e^{-tH} \Phi_j (\sqrt{H}) \Big( \phi_j(\sqrt{H}) f \Big).
\]
The spectrum of the operator $e^{-tH} \Phi_j (\sqrt{H})$ 
is localized around a dyadic number, 
and Lemma~\ref{lem:0321-1} with $G(\lambda)=e^{-t\lambda^2}$, $\lambda>0$ \black implies that there exists a positive constant $C$ such that
\[
\|e^{-tH} \Phi_j (\sqrt{H})\|_{L^p \to L^p} \leq 
C e^{-C^{-1}t 2^{2j}}.
\]
Therefore, we obtain by Proposition~\ref{prop:intro1} (iii) \black that 
\[
\begin{split}
\| e^{-tH}  f \|_{B^{s_2}_{p_2,1}(H)} 
\leq 
& 
C \| e^{-tH}  f \|_{B^{s_2 +d(\frac{1}{p_1}-\frac{1}{p_2})}_{p_1,1}(H)} 
\\
\leq 
&C \sum _{j \in \mathbb{Z}} 
 2^{s_2j + d(\frac{1}{p_1}-\frac{1}{p_2})j}e^{-C^{-1}t 2^{2j}} 
   \| \phi_j(\sqrt{H})f \|_{L^{p_1}}
\\
\leq 
&C \sum _{j \in \mathbb{Z}} 
 2^{ (s_2 - s_1) j + d(\frac{1}{p_1}-\frac{1}{p_2})j}e^{-C^{-1}t 2^{2j}} 
   \| f \|_{B^{s_1 }_{p_1,\infty}}, 
\end{split}
\]
and there exists a positive constant $C$ independent of $t$ such that
\[ 
\sum _{j \in \mathbb{Z}} 
 2^{ (s_2 - s_1) j + d(\frac{1}{p_1}-\frac{1}{p_2})j}e^{-C^{-1}t 2^{2j}} 
 \leq C t^{- \frac{s_2-s_1}{2} -\frac{d}{2}(\frac{1}{p_1} - \frac{1}{p_2}) }.
\]
Here, the convergence of the series follows from the positivity of 
$(s_2 - s_1) + d(1/p_1-1/p_2)$.  
\end{proof}

\begin{proof}[Proof of Theorem 1.3 \black] 
The argument is similar to the proof of Theorem~1.2 in \cite{Iw-2018}. 
We provide only a few comments on the proof. \\
(i) Since $1 \leq q < \infty$, it is crucial that any function can be approximated by 
a finite sum of $\phi_j (\sqrt{H}) f$ over a finite subset of $\mathbb{Z}$. 
We then establish the continuity by applying Lemma~\ref{lem:0321-1} 
to the finite sum. A density argument completes the proof of (i). \\
(ii) When $q = \infty$, the continuity in the dual weak sense 
reduces to the case $q = 1$, where the continuity has already been established in (i).
\end{proof}

\begin{proof}[Proof of Theorem 1.4 \black]
The argument is similar to the proof of Theorem~1.3 in \cite{Iw-2018}. 
The starting point is to establish the following inequality. 

Let $\alpha > 0$, $s_0 \in \mathbb{R}$, and $1 \leq p \leq \infty$. 
Then, there exists a constant $C > 0$ such that 
\begin{equation}\notag 
\begin{split}
& 
C^{-1} (t 2^{\alpha j} ) ^{s_0} e^{-C t 2^{\alpha j} }
\big\| \phi_j (\sqrt{H}) f \big\|_{L^p}
\\
& \leq 
\big\| (tH ^{\frac{\alpha}{2}}) ^{s_0} e^{-t H^{\frac{\alpha}{2}} \black} \phi_j (\sqrt{H}) f 
\big\|_{L^p} 
\leq C ( t 2^{\alpha j} ) ^{s_0} e^{-C^{-1} t 2^{\alpha j} }
\big\| \phi_j (\sqrt{H}) f \big\|_{L^p}
\end{split}
\end{equation}
for any $t > 0$, $j \in \mathbb{Z}$, and $f \in L^p (\mathbb R^d)$. 

The above inequality is established in the same manner as Lemma~5.1 
in \cite{Iw-2018}, using Lemma~\ref{lem:0321-1}. 
We then proceed as in the proof presented in Section~5 of \cite{Iw-2018}. 
\end{proof}

\begin{proof}[Proof of Theorem 1.5 \black]
The proof follows the same argument as that of Theorem~1.4 in \cite{Iw-2018} 
(see Section 6). 
\end{proof}


\noindent
{\bf Acknowledgement.} The first author was supported by JSPS KAKENHI  
Grant Numbers 20K03669.  



\begin{thebibliography}{99}

\bibitem{bcd} H. Bahouri, J.-Y. Chemin, and R. Danchin.
\newblock Fourier analysis and nonlinear partial differential equations. Grundlehren Math. Wiss., 
343 [Fundamental Principles of Mathematical Sciences] Springer, Heidelberg, 2011, xvi+523 pp.

\bibitem{dbdf1} A. de Bouard, A. Debussche and R. Fukuizumi, 
\newblock 
{``Long time behavior of Gross-Pitaevskii equation at positive temperature,''}
\newblock
{SIAM. J. Math. Anal.}
{\bf 50} (2018) 5887-592.

\bibitem{dbdf2} A. de Bouard, A. Debussche and R. Fukuizumi, 
\newblock 
{``Two-dimensional Gross-Pitaevskii equation with space-time white noise,''}
\newblock
{Int. Math. Res. Not.} 
(2023) 10556-10614.

\bibitem{dbdf3} 
A. de Bouard, A. Debussche, and R. Fukuizumi,
\newblock ``Stationary martingale solution for the 2d stochastic Gross-Pitaevskii equation,''
\newblock{RIMS K\^{o}ky\^{u}roku Bessatsu B95 (2024), 23-36.} 

\bib{BD-2015}{article}{
   author={Bui, The Anh},
   author={Duong, Xuan Thinh},
   title={Besov and Triebel-Lizorkin spaces associated to Hermite operators},
   journal={J. Fourier Anal. Appl.},
   volume={21},
   date={2015},
   number={2},
}

\bib{BDY-2012}{article}{
   author={Bui, Huy-Qui},
   author={Duong, Xuan Thinh},
   author={Yan, Lixin},
   title={Calder\'{o}n reproducing formulas and new Besov spaces associated with
   operators},
   journal={Adv. Math.},
   volume={229},
   date={2012},
   number={4},
   pages={2449--2502},
}

\bibitem{dp} P. D'ancona and V. Pierfelice, 
\newblock 
{``On the wave equation with a large rough potential,"}
\newblock{J. Funct. Anal.}
{\bf 227} (2005) 30--77.

\bibitem{d} J. Dziuba\'nski
\newblock 
{``Triebel-Lizorkin spaces associated with Laguerre and Hermite expansions,''}
\newblock
{Proceedings of AMS.}
{\bf 125} (1997) no.12, 3547-3554.

\bibitem{dg} J.Dziuba\'nski and P. Glowacki, 
\newblock 
{``Sobolev spacesrelated to Schr\"odinger operators
with polynomial potentials,''}
\newblock
{Math.Z.} 
{\bf 262} (2009) no.4, 881-894. 

\bibitem{FI} R. Farwig and T. Iwabuchi, 
\newblock
{Sobolev spaces on arbitary domains and semigroups generated by 
the fractional Laplacian}
\newblock
{Bull. Sci. Math.}
{\bf 193} (2024) 103440

\bib{Iw-2018}{article}{
   author={Iwabuchi, Tsukasa},
   title={The semigroup generated by the Dirichlet Laplacian of fractional
   order},
   journal={Anal. PDE},
   volume={11},
   date={2018},
   number={3},
   pages={683--703},
}

\bib{IMT-RMI}{article}{
   author={Iwabuchi, Tsukasa},
   author={Matsuyama, Tokio},
   author={Taniguchi, Koichi},
   title={Boundedness of spectral multipliers for Schr\"{o}dinger operators on
   open sets},
   journal={Rev. Mat. Iberoam.},
   volume={34},
   date={2018},
   number={3},
   pages={1277--1322},
}
\bib{IMT-2019}{article}{
   author={Iwabuchi, Tsukasa},
   author={Matsuyama, Tokio},
   author={Taniguchi, Koichi},
   title={Besov spaces on open sets},
   journal={Bull. Sci. Math.},
   volume={152},
   date={2019},
   pages={93--149},
}
\bib{iwabuchi}{article}{
   author={Iwabuchi, Tsukasa},
   author={Matsuyama, Tokio},
   author={Taniguchi, Koichi},
   title={Bilinear estimates in Besov spaces generated by the Dirichlet
   Laplacian},
   journal={J. Math. Anal. Appl.},
   volume={494},
   date={2021},
}



\bibitem{jn} A. Jensen and S. Nakamura, 
\newblock
{``$L^p$-mapping properties of functions of Schr\"odinger 
operators and their applications to scattering theory,''}
\newblock
{J. Math. Soc.Japan.} 
{\bf 47} (1995) 253-273. 


\bib{PX-2008}{article}{
   author={Petrushev, Pencho},
   author={Xu, Yuan},
   title={Decomposition of spaces of distributions induced by Hermite
   expansions},
   journal={J. Fourier Anal. Appl.},
   volume={14},
   date={2008},
   number={3},
   pages={372--414},
}


\bib{RS_1996}{book}{
   author={Runst, Thomas},
   author={Sickel, Winfried},
   title={Sobolev spaces of fractional order, Nemytskij operators, and
   nonlinear partial differential equations},
   series={De Gruyter Series in Nonlinear Analysis and Applications},
   volume={3},
   publisher={Walter de Gruyter \& Co., Berlin},
   date={1996},
}

\bibitem{t} L. Thomann, 
\newblock
{``Random data Cauchy problem for supercritical Schr\"odinger
equations,''}
\newblock
{Ann. IHP-AN}
{\bf 26} (2009) 2385-2402.

\bibitem{taylor} Taylor
\newblock
{``Tools for PDEs,''}
\newblock 
{Pseudodifferential Operators, Paradifferential Operators, and Layer Potentials, Math. Surveys Monogr., vol. 81,
American Mathematical Society, Providence, RI} (2000)

\bib{Y-2004}{article}{
   author={Yajima, Kenji},
   author={Zhang, Guoping},
   title={Local smoothing property and Strichartz inequality for Schr\"{o}dinger
   equations with potentials superquadratic at infinity},
   journal={J. Differential Equations},
   volume={202},
   date={2004},
   number={1},
   pages={81--110},
}


\end{thebibliography}
\end{document}